\documentclass[12pt,leqno]{article}

\usepackage[official]{eurosym}
\usepackage{fullpage,amsmath,theorem}
\usepackage{epic,eepic}
\usepackage{amssymb,amsmath}
\usepackage{cases}
\def\llvdash{{\|\hskip-2pt \raise 3pt\hbox{\vrule
height 0.25pt width 0.4cm}}}

\def\alp{\alpha}

\def\sig{\sigma}

\def\l{{\langle}}
\def\r{{\rangle}}

\def\calP{\mathcal P}

\def\upr{\upharpoonright}

\def\oa{{\overline A^{\,\lower 7pt_{\hbox{$\scriptstyle\bet}}
\hbox{$\scriptstyle 0\tau$}}}}

\def\bet{\beta}

\def\llvdash{{\|\hskip-2pt \raise 3pt\hbox{\vrule height 0.25pt
width 0.4cm}}}
\def\gam{\gamma}
\newtheorem{theorem}{Theorem}[section]
\newtheorem{lemma}[theorem]{Lemma}
\newtheorem{corollary}[theorem]{Corollary}
\newtheorem{proposition}[theorem]{Proposition}
{\theorembodyfont{\rmfamily}
\newtheorem{definition}[theorem]{Definition}}
{\theorembodyfont{\rmfamily}
\newtheorem{remark}[theorem]{Remark}}
{\theorembodyfont{\rmfamily}
\newtheorem{claim}{Claim}}
{\theorembodyfont{\rmfamily}
}
{\theorembodyfont{\rmfamily}
}
 \DeclareMathOperator{\dom}{dom}
\DeclareMathOperator{\rng}{rng} \DeclareMathOperator{\cof}{cof}

\newcommand{\pr}{\medskip\noindent\textit{Proof}. }

\newcommand{\lusim}[1]{\smash{\underset{\raisebox{1.2pt}[0cm][0cm]{$\sim$}}
{{#1}}}}

\def\dom{{\rm dom}}

\def\rng{{\rm rng}}

\def\supp{{\rm supp}}

\def\llvdash{{\|\hskip-2pt \raise 3pt\hbox{\vrule
height 0.25pt width 0.15cm}}}

\def\Vdashbks{\hbox{$\Vdash\!\!\!\!{\raise2pt\hbox
{$\scriptscriptstyle\backslash$}}$}}

\begin{document}

\title{On fresh sets in  iterations of Prikry type forcing notions}
\baselineskip=18pt
\author{Moti Gitik\footnote{ We are grateful to participants of TAU Set Theory seminar and in particular to Menachem Magidor and to Omer Ben-Neria for their comments and remarks. The work of the first author was partially supported by ISF grant No. 882/22.} and Eyal Kaplan}

\date{\today}
\maketitle

\begin{abstract}
We examine the existence (and mostly non-existence) of fresh sets in commonly used iterations of Prikry type forcing notions. Results of \cite{RestElm} are generalized.
As an application, a question of a referee of \cite{Eyal-thesis} is answered. In addition stationary sets preservation is addressed.

\end{abstract}

\section{ Introduction}

We continue here the series of papers \cite{RestElm}, \cite{G-K-nonst},\cite{K-full}, \cite{G-K-easton}
on iterations of Prikry type forcing notions and concentrate on existence or non-existence of fresh subsets.

Our framework throughout the paper will be as follows:

\begin{enumerate}
  \item GCH,
  \item $\kappa$ is an inaccessible  cardinal,
  \item  $\l P_\alpha, \lusim{Q}_\beta\mid \alpha\leq \kappa, \beta<\kappa\r$ be either Magidor (full) or non-stationary or Easton support iterations of Prikry type forcing notions.
  \item For every $p\in  P_\kappa$ and $\gamma<\kappa$, let $p(\gamma)$ denotes the $\gamma-$th coordinate of $p$,\\ i.e. $p=\l \lusim{p}(\gamma)\mid \gamma<\kappa\r$.
   $\supp(p)=\{\gamma<\kappa\mid  p\upr \gamma \Vdash \lusim{p}(\gamma)\not = 0_{\lusim{Q}_\gamma}\}$.
   \item For every $\beta<\kappa$, $\l Q_\beta, \leq^*_{Q_\beta}\r$ is forced to be $|\beta|-$closed.\\ Note that for a singular $\beta$ this implies $|\beta|^+-$closure.
    \item For every $\beta<\kappa$, $|\lusim{Q}_\beta|<\kappa$.
  \item If the Easton support is used, then for every $p \in P_\kappa$ and every inaccessible $\alpha\leq \kappa$, $\alpha>|\supp(p)\cap \alpha|$, provided that for every $\beta<\alpha$, $|P_\beta|<\alpha$.

  \item If the non-stationary support is used, then for every $p \in P_\kappa$ and every inaccessible $\alpha\leq \kappa$, $\supp(p)\cap \alpha$ is non-stationary in $\alpha$, provided that for every $\beta<\alpha$, $|P_\beta|<\alpha$.

\end{enumerate}


\begin{definition}(J. Hamkins)
A set $A$ of ordinals is called \emph{a fresh set} if all its initial segments are in the ground model, but $A$ is not there.

\end{definition}

Hamkins \cite{hamkins2001gap} proved that if a forcing has a gap at $\kappa$,\footnote{ $P$ has a gap at $\kappa$, if $P=R*\lusim{Q}$ such that $|R|<\kappa$ and $\lusim{Q}$ is $\kappa+1-$strategically closed.} then no fresh subsets of $\kappa$ are added.
\\We proved in \cite{RestElm} results on non-existence of fresh sets in  iterations of certain Prikry type forcing notions $P_\kappa$, as above.
Here we would like to generalize this results.
\\In addition,
 results regarding preservation of stationary sets when forcing with each of the supports.

 The paper is organized as follows.
In Section 2 we deal with preservations of stationarity. Later sections address fresh subsets. Section 3 deal with non-stationary support iterations. Sections 4 and 5
address Easton and full supports respectively.

\section{On preservation of stationarity in Prikry type extensions}

Let $\kappa$ be an inaccessible cardinal and let $S$ be a stationary subset of $\kappa$.
Suppose that $\l P_\alpha, \lusim{Q}_\beta\mid \alpha\leq \kappa, \beta<\kappa\r$ is an Easton or a full (Magidor) or a non-stationary support iteration of Prikry type forcing notions as in the introduction. Let $G_\kappa$ be a generic subset of $P_\kappa$.
\\We would like to address here a question \emph{whether $S$ remains stationary in $V[G_\kappa]$}.

\begin{remark}\label{rem101}
	
	\begin{enumerate}
		\item  Note that we can assume that $\kappa$ is a regular cardinal in $V$, since otherwise it is possible to replace it by $\cof(\kappa)$.
		\item Our main interest will be in situations where $\kappa$ is an inaccessible cardinal in $V[G]$.
		\item If $\kappa$ is a Mahlo cardinal, then Easton support iterations satisfy $\kappa-$c.c. and so preserve stationary subsets of $\kappa$.
		\item If  a full support iteration of Prikry forcings is used, then the set of former measurables which changed their cofinality will be non-stationary, as witnessed by the regressive function which maps each one of them to the first element in its Prikry sequence.
		
	\end{enumerate}
	
\end{remark}

\begin{theorem}\label{thm101}
Assume GCH. Let $\kappa$ be an inaccessible cardinal. Let $P = P_{\kappa}$ be as in the introduction.
Let $G\subseteq P_{\kappa}$ be generic over $V$. Then $\kappa$ is inaccessible in $V\left[G\right]$, and--
	\begin{enumerate}
		\item If $S\subseteq \kappa$ is stationary in $V$ and consists of singulars, then $S$ remains stationary in $V\left[G\right]$.
		\item If the Easton or nonstationary support is used, then $P_{\kappa}$ preserves stationary subsets of $\kappa$.
		\item If the full support is used, $S\subseteq \kappa$ is stationary in $V$, and--
		\begin{enumerate}
			\item For every $\alpha\in S$, $\langle \lusim{Q}_{\alpha}, \leq^* \rangle$ is $\left| \alpha \right|^+$-complete.
			\item For every $\alpha<\kappa$, $\lusim{Q}_{\alpha}$ has the property that for every $p,q,r\in Q_{\alpha}$, if  $p,q\geq^* r$, then there is $t\in Q_{\alpha}$ such that $t\geq^* p,q$.
		\end{enumerate}
		Then $S$ remains stationary in $V[G]$.
	\end{enumerate}

\end{theorem}

The inaccessibility of $\kappa$ in $V\left[G\right]$ is proved in \ref{Cor2-4}. Point 1 of the Theorem is proved in   \ref{Theorem101firstpoint}. Point 2 of the Theorem is proved in Theorems  \ref{Theorem101secondpoint} and \ref{thm ns-new}. Point 3 of the Theorem is proved in \ref{Theorem101fullsupp}.

Let us start with the following lemma which is a weak form of a strong Prikry condition:

\begin{lemma}\label{lem101}
	Let $D\subseteq P_\kappa$ be a dense open and let $p\in P_\kappa$.
	Then there are $\alpha<\kappa$ and $q\geq^* p$ such that for every $r\in P_\alpha, r\geq_{P_\alpha} q\upr \alpha$ there is $r'\geq_{P_\alpha} r$ such that
	$r'{}^\frown q\setminus \alpha \in D$.
	
\end{lemma}
\pr Let $D$ and $p=\l p_\gamma\mid \gamma<\kappa\r$ be as in the statement of the lemma.
\\We prove the lemma for the full support iteration. The arguments for the Easton and for the non-stationary support iterations are very similar, only coordinates in supports should be considered.

Suppose that the conclusion of the lemma fails.
\\We define by recursion, for each $\beta<\kappa$, a condition--
$$p(\beta)=\l \lusim{p}^*_\gamma \mid \gamma<\beta \r^\frown \l \lusim{p}_\gamma \mid \beta\leq\gamma<\kappa\r$$
so that $p(\beta)\upr \beta \Vdash \neg \sigma_\beta$, where
$$\sigma_\beta\equiv  \exists t\in P_\kappa\setminus \beta, t\geq^*p\setminus \beta \exists r\in \lusim{G}_\beta \quad r^\frown t \in D.$$
For the first stage, we have--
$$\sigma_1\equiv  \exists t\in P_\kappa\setminus 1, t\geq^*p\setminus 1 \exists r\in \lusim{G}_1 \quad r^\frown t \in D.$$
$Q_0$ satisfies the Prikry condition, so there is $p_0^*\geq^* p_0$ which decides $\sigma_1$. If $p_0^*\Vdash \sigma_1$, then  $p_0^*{}^\frown \lusim{t}$ will be as desired.
So, assume that $p_0^*\Vdash \neg\sigma_1$, and continue.

The successor step is similar to the first stage above.

For limit steps, suppose that $\beta$ is a limit ordinal. Let us show that--
$$p(\beta)=\l \lusim{p}^*_\gamma \mid \gamma<\beta \r^\frown \l \lusim{p}_\gamma \mid \beta\leq\gamma<\kappa\r$$
is as desired, i.e. $p(\beta)\upr \beta \Vdash \neg \sigma_\beta$.
Suppose otherwise,  then there is $r=\l
\lusim{r}_\gamma\mid\gamma
<\bet\r\in P_\bet$ such that $r\ge p(\bet)\upr\bet$
and $r\Vdash\sig_\bet$.  Extend it, if necessary, so that for some $\lusim{t}$

$$r \Vdash \lusim{t}\geq^*p\setminus \beta  \text{ and } r^\frown \lusim{t} \in D.$$

By the definition of order on
$\calP_\bet$, there is $\bet^*<\bet$  such that
for every $\gam$, $\bet^*\le\gam <\bet$,

$$r\!\!\upr\!\!\gam\Vdash
\lusim{r}_\gam\ge^*_\gam
\lusim{p}^*_\gam.$$  Consider a $P_{\bet^*}$-name
$$\lusim{t}'=\l \lusim{r}_\gamma \mid \beta^*\leq \gamma<\beta\r^\frown \lusim{t}.$$
Then
$$r\upr \beta^*\Vdash \lusim{t'}\geq^* p\setminus \beta^* \text{ and } r\upr \beta^*{}^\frown \lusim{t}' \in D.$$
But $r\upr \beta^*\geq p(\beta^*)\upr \beta^*\Vdash \neg \sigma_{\beta^*}$.
Contradiction.

This completes the construction. \\Consider
$p(\alp)=\l \lusim{p}^*_\gam
\mid\gam<\kappa\r$.  Pick some $r\ge p(\kappa)$ in $D$.  Now we obtain a contradiction as in the limit stage above.
\\
$\square$

Assuming that $\kappa$ is measurable in the ground model, a stronger version of Lemma \ref{lem101} can be proved:

\begin{lemma}\label{lem101withameasurable}
	Assume that $\kappa$ is measurable in $V$ and $U$ is a normal measure on $\kappa$. Let $D\subseteq P_\kappa$ be a dense open and let $p\in P_\kappa$.
	Then there are $\alpha<\kappa$ and $q\geq^* p$, such that $q\upr \alpha=p\upr \alpha$, and for every $r\in P_\alpha, r\geq_{P_\alpha} p\upr \alpha$ there is $r'\geq_{P_\alpha} r$ such that
	$r'{}^\frown q\setminus \alpha \in D$.
	Moreover, for every $X\in U$, $\alpha,q$ above can be chosen such that $\alpha\in X$.	
\end{lemma}

\pr Let $\sigma_{\beta}$ be as in lemma \ref{lem101}. If there exists $\beta\in X$ such that $p\restriction \beta \Vdash \sigma_{\beta}$, we are done.

Assume otherwise.
For every $\beta\in X$ there exists $r_{\beta} \geq p\restriction \beta$ such that $r_{\beta}\Vdash \neg \sigma_{\beta}$. For each such limit $\beta$, there exists $\beta'< \beta$ such that $r_{\beta}\restriction \beta' \Vdash r_{\beta}\setminus \beta' \geq^* p\restriction \left[ \beta', \beta \right)$. The function $\beta\mapsto \beta'$ is regressive, and thus there exist a set $A\in U$ and $\beta^* < \kappa$ such that for every $\beta\in A$, $r_{\beta}\restriction \beta^* \Vdash r_{\beta}\setminus \beta^* \geq^* p\restriction \left[ \beta^*, \beta \right)$. Since $\left| P_{\beta^*} \right|<\kappa$, we can shrink $A\in U$ further, and assume that there exists $r^*\in P_{\beta^*}$, such that, for every $\beta\in A$, $r_{\beta}\restriction \beta^* = r^*$.

Then $r^*$ has the following property: for every $\beta\in A$ there exists $s(\beta)\geq^* p\restriction \left[ \beta^*,\beta  \right)$ such that ${r^*}^{\frown} s(\beta) \Vdash \neg \sigma_{\beta} $. Now apply ineffability: we can find $A^*\subseteq A$, $A^* \in U$, and a $P_{\beta^*}$-name for a condition $s^* = \left[ \beta\mapsto s(\beta) \right]_U \in P\setminus \beta^* $, such that $r^*\Vdash s^*\geq^* p\setminus \beta^*$, and, for every $\beta\in A^*$,
$$ {r^*}^{\frown} s^*\restriction \left[ \beta^*, \beta \right)  = {r^*}^{\frown} s(\beta) \Vdash \neg \sigma_{\beta}. $$
Finally, pick some $q\geq {r^*}^{\frown}s^*$, $q\in D$. Let $\beta\in A^*\setminus \beta^*+1$ be such that-- $$ q\restriction \beta \Vdash q\setminus \beta \geq^* s^*\setminus \beta $$
and in particular,
$$ q\restriction \beta \Vdash q\setminus \beta \geq^* p\setminus \beta, \mbox{ and } {q\restriction \beta}^{\frown} q\setminus \beta \in D $$
and thus $q\restriction \beta \Vdash \sigma_{\beta}$; however, since $\beta\in A^*$,
$$q\restriction \beta \geq {r^*}^{\frown} s^*\restriction{\left[ \beta^*, \beta \right)} \Vdash \neg \sigma_{\beta} $$
which is a contradiction.\\
$\square$

\begin{lemma} \label{Lem2-3}
 Let $p\in P_{\kappa}$ be a condition, and assume that $\lusim{\zeta}$ a $P_{\kappa}$-name for an ordinal. Then there exists $q\geq^* p$ and a set of ordinals $A\in V$ with $\left| A \right|<\kappa$ such that $q\Vdash \lusim{\zeta} \in A$.
\end{lemma}

\pr
Apply lemma \ref{lem101} on the dense open set $D$ of conditions which decide the value of $\lusim{\zeta}$. Then there exists $q\geq^* p$ and $\alpha<\kappa$ such that, for every $r\geq q\restriction \alpha$ there exists $r'\geq r$ such that ${r'}^{\frown} q\setminus \alpha \parallel \lusim{\zeta} $. Let $A\in V$ be the set of all possible values of $\lusim{\zeta}$ as decided by some extension of $q$. We argue that $\left| A\right|<\kappa$.

Assume that $q'\geq q$ decides the value of $\lusim{\zeta}$. Denote $r = q'\restriction \alpha$. Then there exists $r'\geq r$ such that ${r'}^{\frown} q\setminus \alpha \parallel \lusim{\zeta}$. So both   conditions $q'$, ${r'}^{\frown} q\setminus \alpha $ decide the value of $\lusim{\zeta}$; but those conditions are compatible, since $r'\geq p'\restriction \alpha$, and $q'\setminus \alpha \geq q\setminus \alpha$. This shows that every element of $A$ can be realized as the decided value of $\lusim{\zeta}$ by a condition of the form $r^{\frown} q\setminus \alpha$ for some $r\in P_{\alpha}$. But the cardinality of the set of such conditions is strictly below $\kappa$, since  $\left| P_{\alpha} \right| < \kappa$, by the assumption on the cardinality of the forcings $\lusim{Q}_{\beta}$ for $\beta<\kappa$.
\\$\square$

\begin{corollary} \label{Cor2-4}	
Let $G_\kappa\subseteq P_\kappa$. Then $\kappa$ remains inaccessible in $V[G_\kappa]$.
\end{corollary}

\pr
We concentrate on the proof that $\kappa$ remains a regular cardinal after forcing with $P_{\kappa}$, since it's routine to verify that it remains strong limit.

Assume that $\lusim{f}$ is a $P_{\kappa}$-name for a function from some ordinal $\tau<\kappa$ to $\kappa$. Let $p\in P_{\kappa}$ be a condition which forces this. We argue that there exists $p^* \geq p$ and some $\mu^*< \kappa$, such that $p^*\Vdash \rng( \lusim{f} )\subseteq \mu^*$.

Let $G_{\tau+1}\subseteq P_{\tau+1}$ be an arbitrary generic extension containing $p\restriction \tau+1$. We prove that, in $V\left[ G_{\tau+1} \right]$, there exist $q \geq^* p\setminus \tau+1$ and $\mu<\kappa$ such that $q\Vdash \rng(\lusim{f})\subseteq \mu$. Once we prove that, we are done: let $\lusim{q}$, $\lusim{\mu}$ be $P_{\tau+1}$-names which are forced by $p\restriction \tau+1$ to have the above properties. Let $\mu^*<\kappa$ be an upper bound on the set of possible values of $\lusim{\mu}$, as forced by extensions of $p\restriction \tau+1$. Since $\left| P_{\tau+1} \right|<\kappa$, this set is bounded in $\kappa$, and thus there exists such an upper bound below $\kappa$.  Then  $p^* = {p\restriction \tau+1}^{\frown} \lusim{q} \Vdash \rng( \lusim{f})\subseteq \mu^* $, as desired.

Work in $V\left[G_{\tau+1}\right]$. Apply lemma \ref{Lem2-3} over and over to construct a $\leq^*$-increasing sequence of conditions $\langle p_{\xi} \colon \xi \leq\tau \rangle$ in $P_{\kappa}\setminus \tau+1$, such that, for each $\xi<\tau$ there exists some $\mu_{\xi}<\kappa$ such that $p_{\xi} \Vdash \lusim{f}(\xi) < \mu_{\xi}$. Note that in limit steps (including the last step) we may take upper bound, since the direct extension order of $P_{\kappa}\setminus \tau+1$ is more than $\tau$-closed. Finally, $q = p_{\tau}$ forces that the image of $\lusim{f}$ is bounded by $\mu = \bigcup_{\xi<\tau} \mu_{\xi}<\kappa$.
\\$\square$


\begin{theorem} \label{Theorem101firstpoint}
Let $S\subseteq \kappa$ be a stationary set consisting of singulars. Let $G_\kappa\subseteq P_\kappa$. Then $S$ remains stationary in $V[G_\kappa]$.
\end{theorem}

\pr
Let $C\subseteq \kappa$ be a club in $V\left[G\right]$. Let $p\in G$ be condition which forces this.
\\Work in $V$. Pick an elementary submodel $M\preceq H_\chi$ such that:

\begin{enumerate}
	\item $|M|=\delta<\kappa$,
	\item $M\cap \kappa=\delta$,
	\item $\delta\in S$ (in particular, $\delta$ is singular),
	\item ${}^{\cof(\delta)>}M\subseteq M$,
	\item $\kappa,P_\kappa, S, \lusim{C}, p\in M$.
\end{enumerate}
Pick a cofinal sequence in $\delta$, $\l \delta_i \mid i<\cof(\delta) \r$.

Apply lemma \ref{Lem2-3}. Construct (in $V$) a $\leq^*$-increasing sequence of conditions in the forcing $P_{\kappa}\setminus \cof(\delta)$, $\langle p_{\xi} \colon \xi \leq \cof(\delta) \rangle$, such that each condition $p_{\xi}$ belongs to $M$.

We first construct $p_0 \geq^* p$ in $M$, such that $p_0 \restriction \cof(\delta)+1 = p\restriction \cof(\delta)+1$, and, using Lemma \ref{Lem2-3}, \\$p_0 \restriction \cof(\delta)+1 \Vdash \exists \alpha<\kappa, \ p_0 \setminus \cof{\delta}+1 \Vdash \min\left( \lusim{C}\setminus \delta_0 \right) < \alpha$.\\ Let $\alpha_0$ be the least upper bound of the set of all possible values for the $P_{\cof(\delta)+1}$-name $\lusim{\alpha}$. Then $p_0 \Vdash \min\left( \lusim{C}\setminus \delta_0 \right) < \alpha_0$. Note that by elementarity, $\alpha_0 < \delta$.

Assuming that $i<\cof(\delta)$ and $p_{i}, \alpha_i$ have been defined and both are in $M$, and let us define $p_{i+1}, \alpha_{i+1}$. Let $p_{i+1}\geq^* p_{i}$ in $M$ be such that $p_{i+1}\restriction \cof(\delta)+1 = p\restriction \cof(\delta)+1$, and there exists $\alpha_{i+1}< \kappa$,  $p_{i+1} \Vdash \min\left( \lusim{C}\setminus \max\{ \delta_{i+1}, \alpha_i \}  \right) < \alpha_{i+1}$. Take $\alpha_{i+1}$ to be minimal with this property. Then $\alpha_{i+1}< \delta$.

For the limit step, assume that $j< \cof(\delta)$ is limit and $\langle p_{i} \colon i< j \rangle$, $\langle \alpha_i \colon i<j \rangle$ have been constructed. Use the fact that $M$ is closed under $<\cof{\delta}$-sequences to find an upper bound $q\in M$ of $\langle p_i \colon i<j \rangle$, such that $q\restriction \cof(\delta)+1 = p \restriction \cof(\delta)+1 $. We used here the fact that $P_{\kappa}\setminus \cof(\delta)+1 $ is $\mbox{cof}(\delta)$-closed. Finally, let $p_{j}\geq^* q$ be chosen in $M$ such that $p_{j}\restriction \cof(\delta)+1 = p \restriction \cof(\delta)+1 $ and, for some $\alpha_j<\kappa$,  $p_{j}\Vdash \min\left( \lusim{C} \setminus \max\{ \delta_j, \sup\{ \alpha_i \colon i<j \}  \} \right) < \alpha_j$. We used here again the fact that $M$ is closed under $<\cof(\delta)$-sequences and thus $\langle \alpha_i \colon i<j \rangle\in N$. Let $\alpha_j$ be the minimal with the above property. Then $\alpha_j \in M$.

This concludes the inductive construction. In the final limit step, take, in $V$, an upper bound $q^*$ for all the conditions $\l p_{i} \colon i < \cof(\delta) \r$. Then $q^* \Vdash \delta = \bigcup_{i<\cof(\delta) }{ \alpha_i } \in \lusim{C}$, since $C$ is forced by $p$ to be closed. Thus $\delta\in C\cap S$, as desired.
\\$\square$

\begin{remark}
Basically the same argument works for $\kappa$ replaced by $\kappa^+$.

\end{remark}

Now let us try to extend the theorem to $S$'s which consists of regular cardinals. We deal with three supports separately.\\First point out the following:

\begin{theorem} \label{Theorem101secondpoint}
Suppose that  the Easton support is used in $P_\kappa$.
Assume that $S\subseteq \kappa$ is stationary. Then $S$ remains stationary in $V\left[G\right]$.
\end{theorem}

\pr
If $\kappa$ is a Mahlo cardinal, the Easton support iteration is $\kappa$-c.c. and thus preserves stationary subsets of $\kappa$. Thus, we can assume that $\kappa$ is not Mahlo. In this case, there exists a club $C\subseteq \kappa$ of singular cardinals. Thus, by shrinking $S$ we can assume that it consists of singulars.
Then, Theorem \ref{Theorem101firstpoint} applies.\\$\square$

Turn to the full support.

\begin{theorem} \label{Theorem101fullsupp}
Suppose that $P_\kappa$ is   the full support iteration. \\Assume that for every $\beta<\kappa$, for every $p,q,r\in \lusim{Q}_\alpha$, if $p,q\geq^*r$, then there is $t \in \lusim{Q}_\alpha$ such that $t\geq^*p,q$.\\
Let $S\subseteq \kappa$ be a stationary such that,
for every $\alpha\in S$, $\lusim{Q}_\alpha$ is $|\alpha|^+-$complete. \\Then $S$ remains stationary in $V\left[G\right]$.
\end{theorem}

\pr
Suppose otherwise. Pick some $p\in P_\kappa$ and a name $\lusim{C}$ such that

$$p\Vdash \lusim{C} \text{ is a club in } \kappa \text{ and } \lusim{C}\cap S=\emptyset.$$
Pick now $M\preceq H_\chi$ and $\delta\in S$, as in Theorem \ref{Theorem101firstpoint}.
\\It is enough to find a condition $r\geq p$ which forces $``\lusim{C} \text{ is unbounded in } \delta$'' in order to derive a contradiction.
\\Suppose that there is no such $r$. Let $G_\delta\subseteq P_\delta$ be a generic with $p\upr \delta\in G_\delta$.
Then , in $V[G_\delta]$, $$p\setminus \delta\Vdash \lusim{C}_{G_\delta}\cap \delta \text{ is bounded in } \delta.$$
We have $\delta\in S$, so, by the assumption of the theorem, $\lusim{Q}_\delta$ is $\delta^+-$complete, and hence, $P_\kappa/G_\delta$ is $\delta^+-$complete.
Then there are $p'\in P_\kappa/G_\delta, p'\geq^* p\setminus \delta$ and $\rho<\delta$ such that
$$p'\Vdash \lusim{C}_{G_\delta}\cap \delta\subseteq \rho.$$
Pick some $t \in G_\delta, t\geq p\upr \delta$ such that
$$t^\frown \lusim{p}'\Vdash  \lusim{C}_{G_\delta}\cap \delta\subseteq \rho.$$
Since $t\geq p\upr \delta$, there exists $\gamma<\delta$ such that  $t=t\upr \gamma^\frown t\setminus \gamma$ and $t\setminus \gamma\geq^*p\upr [\gamma,\delta)$. In particular, $t\restriction \gamma \Vdash {t\setminus \gamma}^{\frown} \lusim{p}'\geq^* p\setminus \gamma$.
\\Now work in $M$ above the condition $p$. Note that $\gamma, \rho$ above are below $\delta$ and thus are in $M$. Also $p,\lusim{C}\in M$, and recall that
$$p\Vdash \lusim{C} \mbox{ is a club in } \kappa$$
Let $\lusim{\zeta}\in M$ be a $P_\kappa$-name such that $p\Vdash \lusim{\zeta} = \min( \lusim{C}\setminus \check{\rho}+1 ) $. By Lemma \ref{Lem2-3}, we can find (in $M$) a $P_{\gamma}$-name for a condition $\lusim{q}_1\geq^* p\setminus \gamma$ such that $p\upr \gamma \Vdash \exists \mu<\kappa \   \lusim{q}_1\Vdash \lusim{\zeta}< \mu$. Let $\lusim{\mu}\in M$ be a $P_{\gamma}$-name such that $p\upr \gamma \Vdash    \lusim{q}_1\Vdash \lusim{\zeta}< \lusim{\mu}$. In $M$, let $\mu^*$ be the supremum of all possible values of $\lusim{\mu}$, as forced by extensions of $p\upr \gamma$. Then  $\mu^*<\delta$ (since $\mu^* \in M\cap \kappa$). By elementarity, ${p\upr \gamma }^{\frown} q_1 \Vdash \lusim{\zeta}< \mu^*$ holds in $V$ as well.

Finally, let $p^*\in P_{\kappa}$ be a condition such that $q^*\upr \gamma = t\upr \gamma$ and $p^*\setminus \gamma$ direct extends both the conditions ${t\setminus \gamma}^{\frown}\lusim{p}', \lusim{q}_1$. Note that, since ${t\setminus \gamma}^{\frown}\lusim{p}'$ and $\lusim{q}_1$ direct extend $p\setminus \gamma$, the extra assumption of the theorem allows us to construct such $p^*$.

Then on the one hand, $p^*$ forces that $\lusim{C}\cap \delta \subseteq \rho$; on the other hand, it forces that $\min\left( \lusim{C}\setminus \rho+1 \right) \leq \mu^* < \delta $. A contradiction.
\\
$\square$

Turn now to the non-stationary support.

\begin{theorem}\label{thm ns-new}
Assume that the non-stationary support is used. Let $S\subseteq \kappa$ be stationary. Then $S$ is stationary in $V[G_\kappa]$.

\end{theorem}
\pr Let $\lusim{C}$ be a $P_\kappa-$name and $p\in P_\kappa$, $p\Vdash \lusim{C} \text{ is a club disjoint to } S$.

We construct:
\begin{enumerate}
	\item $\leq^*$-increasing sequence of conditions $\langle p_{i} \colon i<\kappa \rangle $,
	\item  decreasing sequence of clubs $\langle C_i \colon i<\kappa \rangle$, each  is disjoint to the support of $p_i$,
	\item increasing  continuous sequences of ordinals $\langle \nu_i \colon i<\kappa \rangle$, $\langle \alpha_i \colon i<\kappa \rangle$.
\end{enumerate}

The sequence of conditions will be a fusion sequence, in the sense that, for every $i<j$, $p_{i}\leq^* p_j$ and $p_{j}\restriction \nu_i = p_i \restriction \nu_i$.

Let $p_0 = p$ and $\alpha_0=0$. Pick a club $C_0$ disjoint to $\supp(p_0)$. Let $\nu_0=\min(C_0)$.
\\Let $G_{\nu_0}\subseteq P_{\nu_0}$ be a generic with $p_0\upr \nu_0\in   G_{\nu_0}$. Apply Lemma \ref{Lem2-3}, inside $V[G_{\nu_0}]$, and find $p'_1 \in P_\kappa/G_{\nu_0}, p_1'\geq^*p_0\setminus \nu_0$ and $A_1'\subseteq \kappa, |A_1'|<\kappa$ such that $p_1'\Vdash \min(\lusim{C})\in A_1'$.
\\Now back to $V$, we have $|P_{\nu_0}|<\kappa$, hence there is $A_1, |A_1|<\kappa$ such that\\ $p_1=p_0\upr \nu_0^\frown \lusim{p}_1'\Vdash  \min(\lusim{C})\in A_1$.
\\Set $\alpha_1=\sup(A_1)$.

Next, we pick a club $C_1\subseteq C_0$ disjoint to $\supp(p_1)$. Let $\nu_1=\min(C_1)\setminus \max(\nu_0, \alpha_1)$.
\\Let $G_{\nu_1}\subseteq P_{\nu_1}$ be a generic with $p_1\upr \nu_1\in   G_{\nu_1}$. Apply Lemma \ref{Lem2-3}, inside $V[G_{\nu_1}]$, and find $p'_2 \in P_\kappa/G_{\nu_1}, p_2'\geq^*p_1\setminus \nu_1$ and $A_2'\subseteq \kappa, |A_2'|<\kappa$ such that $p_2'\Vdash \min(\lusim{C}\setminus \nu_1+1)\in A_1'$.
\\Now back to $V$, we have $|P_{\nu_1}|<\kappa$, hence there is $A_2, |A_2|<\kappa$ such that\\ $p_2=p_1\upr \nu_1{}^\frown \lusim{p}_2'\Vdash  \min(\lusim{C}\setminus \nu_1+1)\in A_2$.
\\Set $\alpha_2=\sup(A_2)$.

We do the same at each successor stage $i<\kappa$ of the construction.
\\Suppose now that $i<\kappa$ is limit. Set $\alpha_i=\nu_i=\bigcup_{j<i}\nu_j=\bigcup_{j<i}\alpha_j$ and let $p_i$ be the coordinateswise   union of $p_j,j<i$.
Let $C_i=\bigcap_{j<i}C_j$. Pick $\nu_{i+1}=\min(C_i)\setminus \nu_i+1$. Continue as above.

Finally, let $p_\kappa$ be the coordinateswise   union of $p_i,i<\kappa$.
It is in $P_\kappa$ since for every $i<j$, $p_{i}\leq^* p_j$ and $p_{j}\restriction \nu_i = p_i \restriction \nu_i$ and $\Delta_{i<\kappa}C_i$ is disjoint to its support.
\\Now, $$p_\kappa\Vdash \{\alpha_i \mid i \text{ limit }\}\subseteq \lusim{C}.$$
But $ \{\alpha_i \mid i \text{ limit }\}$ is club in $V$, and so, $S\cap  \{\alpha_i \mid i \text{ limit }\}\not= \emptyset$. Contradiction.
\\
$\square$

A similar argument can be using elementary submodel shows preservation of stationary subsets of $\kappa^+$ with non-stationary support iterations.

Note $|P_\kappa|=\kappa$, in the Easton support case, and so, $\kappa^+$ is preserved. By \cite{ben2017homogeneous},  $\kappa^+$ is preserved in non-stationary support iterations.
Let us consider the full support iterations.
Note that the Magidor iteration of Prikry forcing notions \cite{Mag} satisfies $\kappa^+-$c.c.. A bit more general:

\begin{proposition}
Suppose that
for every $\beta<\kappa$, for every $s,t,r \in \lusim{Q}_\beta$, if $s\leq_{\lusim{Q}_\beta}^*t,r$, then there is $e \in \lusim{Q}_\beta, e\geq_{\lusim{Q}_\beta}^*t,r$.
Then $P_\kappa$ satisfies $\kappa^+-$c.c.
\end{proposition}

 In general it turns out that $\kappa^+$ may be collapsed with full support iteration.

\begin{proposition}
Suppose that $\kappa$ is  a measurable cardinal.
Let $P_\kappa$ be the full support iteration of  $Col(\alpha,\alpha^+)=\{f\mid f \in {}^\xi\alpha^+ , \xi<\alpha\}$, for every regular $\alpha<\kappa$. Then $\kappa^+$ is collapsed in $V^{P_\kappa}$.
\end{proposition}
\pr
Let $U$ be a normal measure over $\kappa$.
\\We start with the following claim:

\begin{claim}\label{clm2}
Let $p=\l \lusim{p}(\beta) \mid \beta<\kappa\r \in P_\kappa$.
Then there are $A^*\in U$, $\tau^*<\kappa$ and $p^*\geq p$ such that $p^*\upr \beta \Vdash \dom(\lusim{p}^*(\beta))=\tau^*$, for every $\beta\in A^*$.

\end{claim}
\pr
For every regular $\beta<\kappa$ there are $s_\beta\in P_\beta, s_\beta\geq p\upr \beta$ and $\tau_\beta^0<\beta$ such that $s_\beta\Vdash \dom(\lusim{p}(\beta))=\tau_\beta^0$.
\\
 Find $A_0' \in U$ and $p_0\in P_\kappa$ such that $p_0\upr \beta= s_\beta,$ for every $\beta\in A_0'$. For example, take $p_0=[\beta\mapsto s_\beta]_U$.
 \\Then we consider a regressive function $\beta\mapsto \tau_\beta^0$ on $A_0'$. Find $A_0\subseteq A_0', A_0 \in U$ and $\tau^0<\kappa$ such that $\tau_\beta^0=\tau^0,$ for every $\beta\in A_0$.
\\Repeat the process with $p_0$ replacing $p$ and find $A_1\subseteq A_0, A_1 \in U$, $p_1$ and $ \tau^1$
such that $p_1\upr \beta \Vdash \dom(\lusim{p}_0(\beta))=\tau^1$.
\\Continue by induction. Let $A^*=\bigcap_{n<\omega}A_n$ and $p^*$ be the coordinatewise  union of $p_n$'s.
Set $\tau^*=\bigcup_{n<\omega}\tau^n$.
\\Then $p^*\upr \beta \Vdash \dom(\lusim{p}^*(\beta))=\tau^*$, for every $\beta\in A^*$, will be as desired.
\\
$\square$ of the claim.

For every $\tau<\kappa$, define a maximal antichain $A_\tau$  in $P_\kappa$.
Proceed as follows.
\\Let us pick functions $\l h_\gamma\mid \gamma<\kappa^+\r$ such that $\dom(h_\gamma)\in U$, for every $\alpha\in \dom(h_\gamma), h_\gamma(\alpha)<\alpha^+$ and $[h_\gamma]_U=\gamma$, for example $\kappa^+-$canonical functions will do the job.
\\Fix $\tau<\kappa$.
Let $A_\tau$ be a maximal antichain in $P_\kappa$ of cardinality $\kappa^+$ which consists of $p\in P_\kappa$ such that:

(*)\emph{if for some $B\in U$ and $\gamma<\kappa^+$,  the condition $t_{B \gamma}=\l {t}_{B \gamma}(\alpha) \mid \alpha<\kappa\r$ is compatible with $p$, then, for some $B'\in U$, $p\geq t_{B' \gamma}$,}
\\where, for $E\in U$, ${t}_{E \gamma}(\alpha)=0_\alpha,$ unless $\alpha\in E\cap \dom(h_\gamma)\setminus \tau+1$, and if  $\alpha\in E\cap \dom(h_\gamma)\setminus \tau+1$, then
${t}_{E\gamma}(\alpha) =\{(\tau, h_\gamma(\alpha))\}$, i.e. the value of the generic function for $\alpha$ at $\tau$ is $h_\gamma(\alpha)$.

Let $\l p^\tau_i\mid i<\kappa^+\r$ be an enumeration of $A_\tau$.

Let $G\subseteq P_\kappa$ be a generic. Define $F:\kappa\to (\kappa^+)^V$ by setting $F(\tau)=i$ iff $p^\tau_i\in G$.
\\We claim that such $F$ is onto. \\Suppose otherwise.
Pick some $p\in G$ and $\eta<\kappa^+$ such that $p\Vdash \rng(\lusim{F})\subseteq \eta$.\\
Apply Claim \ref{clm2}. Let $A^*\in U$, $\tau^*<\kappa$ and $p^*\geq p$ be as in the conclusion of the claim.
\\Now for every $\gamma<\kappa^+$, we can extend $p^*$ to a condition $p^\gamma$ by adding a pair $(\tau^*, h_\gamma(\alpha))$ to $p^*(\alpha)$, for every $\alpha\in A^*\cap \dom(h_\gamma)\setminus \gamma+1$.
 Note that if $\gamma\not = \gamma'$, then $p^\gamma,p^{\gamma'}$ are incompatible. So, the set $\{ p^\gamma\mid \gamma<\kappa^+\}$
consists of $\kappa^+-$many incompatible conditions.
\\Then, each of $p^\gamma$'s must be compatible with a member of $A_{\tau^*}$ with index below $\eta$.
Hence, there is $i^*<\eta$ such that $p^{\tau^*}_{i^*}$ is compatible with $\kappa^+-$many $p^\gamma$'s.
Pick two of them $\gamma\not = \gamma'$.
\\By (*), then  $p^{\tau^*}_{i^*}\geq t_{B \gamma}, t_{B' \gamma'}$, for some $B,B'\in U$.
However, $\gamma\not = \gamma'$ implies that there is  $\alpha\in B\cap B'$, such that $h_\gamma(\alpha)\not =h_{\gamma'}(\alpha)$, and this is impossible due
 to the compatibility.
\\Contradiction.
\\
$\square$.

We conjecture that the measurability assumption can be much weaken. However, the following positive result can be proved:

\begin{proposition}
Suppose that there is a club $C\subseteq \kappa$ such that for every $\alpha\in C$, $\l \lusim{Q}_\alpha,\leq^*_{\lusim{Q}_\alpha}\r$ is forced to be $\alpha^+-$closed.
Then $P_\kappa$ preserves $\kappa^+$.
\\In particular, if $\kappa$ is not a Mahlo cardinal, then $P_\kappa$ preserves $\kappa^+$.
\end{proposition}
\pr
 Let $p \in P_\kappa$ and $\lusim{f}$ be a name such that

$$p\Vdash \lusim{f}:\kappa \to \kappa^+.$$

Fix a club $C$ such that for every $\alpha\in C$, $\l \lusim{Q}_\alpha,\leq^*_{\lusim{Q}_\alpha}\r$ is forced to be $\alpha^+-$closed. Assume also that for every $\alpha\in C, $ for every $\beta<\alpha, |P_\beta|<\alpha$.
Let $\l \alpha_i \mid i<\kappa\r$ be an increasing continuous enumeration of $C$.
Apply Lemma \ref{Lem2-3} and find $p_0\geq^* p, p_0\upr \alpha_0=p\upr \alpha_0$ and $\eta_0<\kappa^+$ such that $p_0\Vdash \lusim{f}(0)<\eta_0$.
\\Continue by induction and define a $\leq^*-$increasing sequence $\l p_i \mid i<\kappa \r$ and sequence $\l \eta_i \mid i<\kappa \r$ of ordinals below $\kappa^+$ such that
\begin{enumerate}
  \item $p_i\Vdash \lusim{f}(i)<\eta_i$,
  \item for every $i<j$, $p_j\upr \alpha_i=p_i\upr \alpha_i$,

\end{enumerate}

There is no problem at limit stages $i$, since $\l P_\kappa\setminus \alpha_i, \leq^*\r$ is $\alpha_i^+-$closed since $\alpha_i\in C$.

The second item insures that there is $p^*\in P_\kappa$ such that $p^*\geq^* p_i,$ for every $i<\kappa$.
Then
$$p^*\Vdash \rng(\lusim{f})\subseteq \bigcup_{i<\kappa}\eta_i.$$
\\
$\square$

\section{ Non-stationary support iterations }

We assume GCH throughout as before.
Let $\kappa$ be an inaccessible cardinal.

Let $\langle  P_{\alpha}, \lusim{Q}_{\beta} \colon \alpha\leq \kappa, \beta<\kappa\rangle$ be a non-stationary support iteration of Prikry type forcings,
with the properties stated in the introduction.

Let $I$ be a stationary subset of $ \kappa$ which consists of singular cardinals $\mu$ such that for every $\gamma<\mu, |P_\gamma|<\mu$.
\\ Assume that- 
\\ \emph{{If $\alpha\in I$, then $\Vdash_{P_{\alpha}} \langle P_\kappa\setminus {\alpha} , \leq^* \rangle \mbox{ is }\alpha^{++} \mbox{-closed.}$}}

Note that for a singular cardinal $\alpha$, we have  $\Vdash_{P_{\alpha}} \langle P_\kappa\setminus {\alpha}, \leq^* \rangle \mbox{ is }\alpha^{+} \mbox{-closed.}$
We will need a slightly more closure.
A typical situation is where $\kappa$ is Mahlo, $Q_\alpha$ is trivial at every accessible $\alpha$ and each forcing $\lusim{Q}_{\alpha}$ has cardinality below the least inaccessible above $\alpha$.

\begin{theorem} \label{thm4-1}
	 $P_\kappa$ does not add fresh unbounded subsets to $\kappa$.
\end{theorem}
\pr
Recall the following  fusion lemma for the non-stationary support iteration of Prikry forcings\footnote{The fusion property for non-stationary support iterations is due to Friedman and Magidor \cite{friedman2009number}. A version suitable for iterations of Prikry type forcings appeared  in \cite{ben2017homogeneous}. The proof is basically given in lemma 3.3 in \cite{RestElm}.}:

\begin{lemma} \label{lem-f1}
	Let $ p\in P_{\kappa} $. For every $ \beta < \kappa $, let $ F(\beta) $ be a $ P_{\beta} $-name for a $ \leq^* $-dense open subset of $ P\setminus \beta $ above $ p\setminus\beta $, and assume that this is forced by $ p\restriction{\beta} $. Then there exist $p^* \geq^{*} p$ and a club $ C\subseteq \kappa $ such that for every singular $\beta \in C$, $ p^* \restriction{\beta} \Vdash p^* \setminus \beta \in F(\beta).$
\end{lemma}

Let $G\subseteq P_{\kappa}$ be generic over $V$, and assume for contradiction that there exists a function $f\in 2^{\kappa}$ which is the characteristic function of a fresh subset of $\kappa$. Let $ \lusim{f} $ be a $ P_{\kappa} $-name for it, and assume that this is forced by some condition in $G$. For simplicity, assume that this is the weakest condition.

Let $\zeta\in \kappa\cap I$ be the least ordinal for which a new subset is added in the extension from $V$ to $V\left[G\right]$. Such $\zeta $ exists, since the forcings $Q_{\alpha}$ (for $\alpha<\kappa$) have cardinality below $\kappa$, and at least one of them is non-trivial.

Note that $\zeta\in \kappa\cap I$,  $\l P_\kappa\setminus \zeta, \leq^* \r$  is  $\zeta^{++}$-closed. Pick a condition $q\in P_{\zeta}$ which forces that a new subset is added to $\zeta$. For simplicity, assume that the weakest condition in $P_{\zeta}$ forces this (else, work above a condition in $P_{\kappa}$ whose restriction to $P_{\zeta}$ equals $q$).

We divide into two cases:

\textbf{Case 1.}\footnote{ It basically repeats those of 4.11, \cite{RestElm}.} There exists $ \mu\in \left( \zeta, \kappa \right)\cap I$ and a condition $ p^* \in P_{\mu} $ which forces that the following property holds:
\begin{align*}
	&\exists p\in G_\mu \  \exists s\in P\setminus \mu \forall r\geq^* s \  \exists \xi <\kappa \  \exists r_0,r_1\geq^* r  , \\
	& V\vDash \left( {p}^{\frown} r_0 \parallel {\lusim{f}}\restriction { \xi }  , \  {p}^{\frown} r_1 \parallel {\lusim{f}}\restriction { \xi } \right) \mbox{, and the decisions are different.}
\end{align*}
(here, $ G_{\mu} $ denotes the canonical name for the generic set for $ P_{\mu} $). By extending $ p^* $, we can decide the value of $ p $ in the statement above, and thus assume that $ p^* \geq p $. Let $ \lusim{s} $ be a $ P_{\mu} $-name for $ s $ from the above property, and assume that this is forced by $ p^* $.


Let us apply the same methods as in the main lemma in \cite{hamkins2001gap}. We construct, in $ V $, a binary tree of conditions, $ \langle  \langle p^*, \lusim{s}_{\sigma} \rangle \colon \sigma\in {}^{\mu>}2\rangle $ and a tree of functions $ \langle b_{\sigma} \colon \sigma\in {}^{\mu>}2 \rangle  $ such that $ \lusim{s}_{\emptyset} = \lusim{s} $, and for every $ \sigma\in {}^{\mu>}2$:
\begin{enumerate}
	\item $\forall i<2, \  \langle p^* , \lusim{s}_{\sigma^{\frown} \langle i \rangle } \rangle \parallel  \lusim{f}\restriction { lh\left( b_{\sigma ^{\frown} \langle i \rangle  } \right) }  = b_{\sigma ^{\frown} \langle i  \rangle} $.
	\item $b_{\sigma^{\frown} \langle 0 \rangle  } \perp b_{\sigma^{\frown} \langle 1 \rangle  }$.
	\item If $ \mbox{lh}(\sigma) $ is limit, then $ p^* $ forces that $ \lusim{s}_{\sigma} $ is an upper bound, with respect to the direct extension order, of $\langle  \lusim{s}_{\sigma\restriction{ \xi}} \colon \xi< \mbox{lh}(\sigma)  \rangle  $.
	\item $ b_{\sigma} $ is an end extension of $ b_{\sigma\restriction { \xi}} $ for every $ \xi <\mbox{lh}(\sigma) $.
\end{enumerate}

Now assume that $ g\subseteq P_{\mu} $ is generic over $ V $ with $ p^* \in g $. In $ V\left[g\right] $, let $ h\in 2^{<\mu} $ be the characteristic function of
a new subset of $\mu$ (such a new subset exists because $\mu$ is above $\zeta$). $ h $ defines a branch through the binary tree, $ \langle \langle p^*, \lusim{s}_{h\restriction { \xi }} \rangle  \colon  \xi < \mu  \rangle  $. Since $\langle \lusim{s}_{h\restriction { \xi }}  \colon \xi<\mu \rangle$ form the sequence of moves of the first player, when the second plays according to the strategy $\tau$, there exists an upper bound $s^*\in P\setminus \mu$, which extends all the conditions in the sequence. Thus, there exists an upper bound for the branch, of the form $ \langle p^*, \lusim{s}^* \rangle  $. It forces that--
$$ b= \bigcup_{\xi<\mu} b_{h\restriction { \xi }}  $$
is an initial segment of $ \lusim{f} $. We argue that this must be a strict initial segment of $f$. Indeed, otherwise,  $\left(\mbox{cf}(\kappa)\right)^{V\left[G_{\mu}\right]} \leq \mu$. But, since $\kappa$ is inaccessible, $G_{\mu}\subseteq P_{\mu}$ is a forcing whose cardinality is strictly below $\kappa$, so it preserves cofinalities greater of equal to $\kappa$.

Therefore, $ b$ is a strict initial segment of $f$, and thus $b\in V$. Therefore, $ h $ can be defined, in $ V $, using the binary tree and the set $ b $. This is a contradiction to the choice of $ h $.

\textbf{Case 2.} For every $ \mu\in \left( \zeta, \kappa \right)\cap I $, every condition in $ P_{\mu} $ forces that--
\begin{align*}
	& \forall p\in \lusim{G}_\mu \ \forall s\in P_\kappa\setminus \mu \exists r\geq^* s \  \forall \xi<\kappa  \  \forall r_0,r_1\geq^* r, \\
	& V\vDash \mbox{ If }  {p}^{\frown} r_0 \parallel {\lusim{f}}\restriction { \xi }  \mbox{ and } \  {p}^{\frown} r_1 \parallel {\lusim{f}}\restriction { \xi }  \mbox{ then the decisions are the same.}
\end{align*}
Define for every  $ \mu<\kappa$ 
\begin{align*}
	e\left(\mu\right) = \{ & r\in P_\kappa\setminus \mu \colon \forall p\in G_{\mu} \     \forall \xi<\kappa \  \forall r_0,r_1\geq^* r, \\
	& V\vDash \mbox{ If }  {p}^{\frown} r_0 \parallel {\lusim{f}}\restriction{ \xi }  \mbox{ and } \  {p}^{\frown} r_1 \parallel {\lusim{f}}\restriction{ \xi }  \mbox{ then the decisions are the same.}   \}
\end{align*}

\begin{claim}
Suppose that $\mu\in \left( \zeta, \kappa \right)\cap I $. Then
$e(\mu)$ is a dense open in $\l  P_\kappa\setminus \mu, \leq^*\r$.\footnote{ It is the only place in the proof where $\mu^{++}-$completeness is used.}

\end{claim}
\pr
Just note  that $ \left| G_{\mu} \right|\leq \mu^+$ and $\l  P_\kappa\setminus \mu, \leq^*\r$ is $\mu^{++}-$closed.
\\
$\square$ of the claim.

Given $\mu\in \left( \zeta, \kappa \right)\cap I$ as above, the following set is also forced to be $ \leq^* $-dense open in $ P\setminus \mu$:
\begin{align*}
	d\left(\mu\right) = \{ & r\in P\setminus \mu \colon \exists g\in 2^{\mu}, \ r\Vdash \lusim{f}\restriction \mu = g\}.
\end{align*}
The $\leq^*$-density of $d(\mu)\subseteq P\setminus \mu$ follows as well from the fact that the direct extension order of $P\setminus \mu$ is more than $\mu$-closed.

We can now apply the standard fusion argument \ref{lem-f1}. There exists $ p\in P_{\kappa} $ and a club $ C\subseteq \kappa $ such that $\min (C)> \zeta $, and, for every $ \mu\in C\cap I$,
$$ p\restriction{ \mu} \Vdash  p\setminus \mu \in  d(\mu) \cap e(\mu). $$

For each $\mu\in C\cap I$,  there exists a condition in $G$ of the form ${q_{\mu}}^{\frown} p\setminus \mu $, where $q_{\mu}\in P_{\mu}$, which decides the value of $\lusim{f}\restriction {\mu}$. The reason is that $p\restriction {\mu}$ forces that $p\setminus \mu$ is in $d(\mu)$, and thus the value of $\lusim{f}\restriction {\mu}$  is decided by the forcing $P_{\mu}$.

For each such $\mu$, $q_{\mu}$ is an extension of $p\restriction{\mu}$, and thus there exists a finite set $b_{\mu}\subseteq \mu$ such that at every $\delta \in \mu \setminus b_{\mu}$, $q_{\mu}(\delta)$ direct extends $p(\delta)$ (as forced by $q_{\mu}\restriction{\delta}$).

The function $\mu\mapsto \max b_{\mu}$ is a regressive function in $V\left[G\right]$, and its domain is the set $C\cap I$. $C\cap I$ is stationary in $V$ since $I$ is assumed to be stationary. By theorem \ref{thm101}, it is also stationary in $V\left[G\right]$. Since $\kappa$ is still regular in $V\left[G\right]$, we can find an unbounded subset $S\subseteq \kappa$ and an ordinal $\mu^*< \kappa$ such that for every $\mu\in S$, $b_{\mu}\subseteq \mu^*$.  By increasing $\mu^*$, we can assume that it belongs to $C\cap I$.

Now, shrink $S$ further to stabilize the function $\mu \mapsto q_{\mu}\restriction{\mu^*}$. This is possible since $S$ is unbounded in $\kappa$, and $q_{\mu}\restriction{\mu^*}$ is a condition in $P_{\mu^*}$ which has a small cardinality (and, again, $\kappa$ is inaccessible in $V\left[G\right]$).

So we can assume that there exists a condition $q^* \in P_{\mu^*}$, such that for every $\mu\in S$, there exists some direct extension $r_{\mu}\in P\setminus \mu^*$ of $p\setminus \mu^*$, such that--
$$ {q^*}^{\frown} r_{\mu} \parallel \lusim{f}\restriction{\mu} $$
and ${q^*}^{\frown} r_{\mu}\in G$.

Recall that $\mu^*\in C\cap I$ and the condition $p$ obtained by fusion above also satisfies that $p\restriction{\mu^*} \Vdash p\setminus \mu^*\in e(\mu^*)$. Since $q^* \geq p\restriction{\mu^*}$, the condition $q^*$ forces that, for every $\xi<\kappa$, any pair of direct extensions of $p\setminus {\mu^*}$ which decide $\lusim{f}\restriction{\xi}$, decide this initial segment the same way.

It follows that ${q^*}^{\frown}  p\setminus \mu^*$ decides $\lusim{f}$ entirely, and forces it to be the following function of $V$:

\begin{align*}
	h = \bigcup_{ \mu \in \kappa\setminus \mu^* } \{  g\in 2^{\mu} \colon &\mbox{there exists a } P_{\mu^*}\mbox{-name for an extension} \lusim{r}\geq^* p\setminus \mu^*, \\
	&\mbox{such that } {q^*}^{\frown}r \Vdash \lusim{f}\restriction{\mu} = g   \}
\end{align*}
which is a contradiction.
$\square$\\

\begin{remark} \label{remark: GCH on kappa in NS}
	Lemma \ref{lem-f1} implies that $V[G]\vDash 2^{\kappa}  =\kappa^{+}$. Indeed, assume that $A$ is a subset of $\kappa$ in $V[G]$. Let $\lusim{A}$ be a  $P_{\kappa}$-name for it. For every singular $\beta<\kappa$, define in $V^{P_{\beta}}$ the set--
	$$ F(\beta) = \{ q\in P_{\kappa}\setminus \beta \colon \exists A_{\beta}\subseteq \beta, \ q\Vdash \lusim{A}\cap \beta = A_{\beta}  \}. $$
	Note that $\beta$ is singular and thus $F(\beta)$ is $\leq^*$-dense open. Let $p^*\in G$ and $C\subseteq \kappa$ be such that for every singular $\beta\in C$, $p^*\upr \beta \Vdash p^*\setminus \beta\in F(\beta)$. Then there exists a $P_{\beta}$-name $\lusim{A}_{\beta}$ for a subset of $\beta$, such that $p^*\Vdash \lusim{A}\cap \beta = (\lusim{A}_{\beta})_{\lusim{G}\upr P_{\beta}}$. Then $A = (\lusim{A})_G$ can be computed in $V[G]$ from the sequence $\l \lusim{A}_{\beta} \colon \beta\in C \r$. By using canonical names for bounded subsets of $\kappa$, and by $\mbox{GCH}$ in $V$, there are at most $\kappa^{+}$ such sequences. So there are at most $\kappa^+-$many subset of $\kappa$ in $V[G]$.
\end{remark}

The situation with higher cardinals was clarified in  \cite{RestElm}.
The following
 was shown basically in \cite{RestElm}, 4.11:

\begin{lemma}\label{lem>}
 $P_{\kappa}$ does not add fresh unbounded subsets to $\kappa^{+}$, or to any cardinal $\lambda$ of $V$ with $\mbox{cf}\left( \lambda \right) >\kappa$.
\end{lemma}

\pr
The proof is a variation of the proof of theorem \ref{thm4-1}, and it basically appears in \cite{RestElm}. Assume that $\lusim{f}$ is a $P_{\kappa}$-name for the characteristic function of a fresh unbounded subset of $\lambda$. Divide into cases as in the proof of theorem \ref{thm4-1}. Case $1$ remains the same. Case $2$ is simplified, since the sets $d(\mu)$ is no longer required. Indeed, in the notations of the proof of theorem \ref{thm4-1}, assume that $p\in P_{\kappa}$ is a condition and $C\subseteq \kappa$ is a club, such that for every $\mu\in C\cap I$,
$$ p\restriction \mu \Vdash p\setminus \mu \in e(\mu) $$
Now, work in $V\left[G\right]$. For each $\xi<\lambda$, let $q_{\xi}\in P_{\kappa}$ be an extension of $p$ which decides $\lusim{f}\restriction \xi$. Let $b_{\xi}\subseteq \kappa$ be a finite set such that, for every $\alpha\in b_{\xi}$, $q_{\xi}\restriction_{\alpha}\Vdash q_{\xi}(\alpha)\geq^* p(\alpha)$. Let $\mu_{\xi}<\kappa$ be an upper bound on $b_{\xi}$. Since $\mbox{cf}(\lambda) >\kappa$ in in $V$, the same holds true in $V\left[G\right]$ as well (by the same proof as in corollary \ref{Cor2-4}). Thus, there exists $\mu^*\in C\cap I$ such that, for an unbounded $S\subseteq \lambda$, $\mu_{\xi} < \mu^*$. By shrinking $S$, we can assume that, for some $q^* \in P_{\mu^*}$, $q_{\xi}\restriction_{\mu^*} = q^*$. Then $q^*$ satisfies that, for every $\xi\in S$, there exists a direct extension $r_{\xi}\in P\setminus \mu^*$ of $p\setminus \mu^*$, such that ${q^*}^{\frown} r_{\xi} \parallel \lusim{f}\restriction \xi $. Now, as in the proof of case $2$ in theorem \ref{thm4-1}, the condition ${q^*}^{\frown}  p\setminus \mu^*$ forces $\lusim{f}$ to be  the following function of $V$:

\begin{align*}
	h = \bigcup_{ \xi < \kappa^{+}} \{  g\in 2^{\mu} \colon &\mbox{there exists a } P_{\mu^*}\mbox{-name for an extension} \lusim{r}\geq^* p\setminus \mu^*, \\
	&\mbox{such that } {q^*}^{\frown}r \Vdash \lusim{f}\upr{\xi} = g   \}
\end{align*}
which is a contradiction.\\
$\square$

We finish with an application for iterations of Prikry type forcings with the nonstationary support. The referee of \cite{Eyal-thesis} asked if such an iteration, below a cardinal $\kappa$, can add new measurable cardinals below $\kappa$. We show  that the answer is negative\footnote{The answer for the same question in the full or Easton support is known to be negative, see \cite{hb}.}.

\begin{theorem} \label{corollary: NS dont add new measurables}
	Assume GCH and let  $\kappa$ be an inaccessible cardinal.
	\\Let $ P_{\kappa}$ be a non-stationary support iteration of Prikry-type forcings  satisfying the conditions from the beginning of the section. Assume:
	
	$(*)$ For every Mahlo cardinal $\alpha<\kappa$ which is not measurable in $V$, $\l P_\kappa \setminus \alpha,\leq^*\r$ is $\alpha^{+}-$closed.\\
	Let  $\lambda$ be a cardinal such that $\forall \tau<\lambda (|P_\tau|<\lambda)$. Let $G\subseteq P_\kappa$ be generic over $V$. Then If $\lambda$ is measurable in $V[G]$, it was already measurable in $V$. \\
	Furthermore, if the assumption $(*)$ is strengthened to--
	
	$(**)$ For every Mahlo cardinal $\alpha<\kappa$, $\l P_\kappa \setminus \alpha,\leq^*\r$ is $\alpha^{+}-$closed.\\
	then $\lambda$ is measurable in $V[G]$ if and only if it is measurable in $V$.
\end{theorem}

\pr
Assume first that $(**)$ holds and $\lambda$ is a measurable cardinal in $V$. If $\lambda> \kappa$, the Levy-Solovay Theorem \cite{Levy-Solovay} shows that $\lambda$ is measurable in $V[G]$. If $\lambda\leq \kappa$, standard arguments show that $\lambda$ remains measurable in $V\left[G\upr P_{\lambda}\right]$ (see \cite{G-K-nonst}). If $\lambda<\kappa$, the forcing $\l P_{\kappa}\setminus \lambda, \leq^* \r$ is $\lambda^+$-closed, so $P_{\kappa}\setminus \lambda$ does not add subsets to $\lambda$ and  $\lambda$ remains measurable in $V[G]$.
Thus, let us concentrate on the other direction assuming $(*)$.

We first recall Lemma $2.1$ in \cite{RestElm}: given a forcing notion which does not add new fresh unbounded subsets to cardinals of $V$ in the interval $\left[ \kappa, \left(2^{\kappa}\right)^{V} \right] = \left[ \kappa, \kappa^{+} \right]$, every $\kappa$-complete ultrafilter in the generic extension extends a $\kappa$-complete ultrafilter from $V$.
Assume now that $\alpha$ is measurable in $V\left[G\right]$, and let $W\in V\left[G\right]$ be a nontrivial $\kappa$-complete ultrafilter on $\alpha$. If $\alpha\geq \kappa$, then by the results in this section, $P_{\kappa}$ does not add fresh unbounded subsets to $\alpha$, $\alpha^{+}$, and thus $W\cap V \in V$ by Lemma $2.1$ in \cite{RestElm}. Thus, assume that $\alpha<\kappa$, and assume  that $\alpha$ is not measurable in $V$.
So, $\alpha$ is a Mahlo cardinal in $V$. If
the forcing $\l P_\kappa \setminus \alpha,\leq^*\r$ is $\alpha^{++}-$closed, then $W\in V\left[G\restriction P_\alpha\right]$.
\\However we assumed only that  $\l P_\kappa \setminus \alpha,\leq^*\r$ is $\alpha^{+}-$closed. In this case $W$ need not be in $V\left[G\restriction P_\alpha\right]$.
\\ Proceed then as follows. Work in $V\left[G\restriction P_\alpha\right]$. Let $\l A_i \mid i<\alpha^+\r$ be an enumeration of all subsets of $\alpha$ (such enumeration exists by applying  remark \ref{remark: GCH on kappa in NS} on $P_{\alpha}$).
Define a $\leq^*-$ increasing sequence of conditions $\l p_i \mid i<\alpha \r$ in  $V\left[G\restriction P_\alpha\right]$ such that for every $i<\alpha$,
$p_i || A_i \in \lusim{W}$. Set $W'=\{A_i \mid i<\alpha, p_i\Vdash A_i \in \lusim{W}\}$.
Then $W'$ will be an $\alpha-$complete ultrafilter over $\alpha$ in $V\left[G\restriction P_\alpha\right]$.

Apply now \ref{thm4-1} and \ref{lem>} to $ P_\alpha$. It follows that no fresh subsets are added to $\alpha, \alpha^+$.
Now, by Lemma $2.1$ from \cite{RestElm}, $W'\cap V\in V$ is a nontrivial $\alpha$-complete ultrafilter over $\alpha$ in $V$. A contradiction.\\
$\square$

\begin{remark}

\begin{enumerate}
  \item The closure assumptions made on $\l P_\kappa \setminus \alpha,\leq^*\r$ are needed.
\\
For example, start with $V=L[U]$, where $U$ is a normal ultrafilter over $\alpha$. Iterate $Cohen(\beta), \beta<\alpha$. Let $V$ be this model.
Then $\alpha$ is not a measurable in $V$. Force with $Cohen(\alpha)$ over $V$. Then $\alpha$ will be a measurable in the extension.
Here we take $Q_\beta$ to be trivial for every $\beta<\alpha$ and $Q_\alpha=Cohen(\alpha)$.

  \item Also, the assumption $\forall \tau<\lambda (|P_\tau|<\lambda)$ is necessary.
Just use the previous example (with $\lambda=\alpha$). $Q_0=Cohen(\alpha)$ resurrects measurability of $\alpha$.

\item Note that a measurable cardinal in $V$ need not be such in $V[G]$ without assuming $(*)$. Just use the Prikry forcing or the iteration of such forcings.

\end{enumerate}

\end{remark}

\section{On fresh sets in the Easton support iterations }

We start with an easier case of the Easton support iterations with $\kappa$ being a Mahlo cardinal.
 $\kappa-$c.c. of the forcing will be used to
 show the following:

\begin{theorem}\label{thm3-1}
Suppose that $\kappa$ is a Mahlo cardinal and $P_\kappa$ is an Easton support iteration.
Then no fresh sets are added to $\kappa$.
\end{theorem}
\pr
Suppose otherwise. Work in $V$. Let $\lusim{A}$ be a name of such subset.

Define a tree $T$ of possibilities as follows.
Fix an increasing enumeration $\l \kappa_\xi\mid \xi<\kappa\r$ of all inaccessible cardinals below $\kappa$.

For every $\xi<\kappa$, let
$$Lev_\xi(T)=\{x\subseteq \kappa_\xi \mid \exists p\in P_\kappa \quad p\Vdash \lusim{A}\cap \kappa_\xi=x\}.$$

Let $x\in Lev_\alpha(T), y \in Lev_\beta(T)$.
Set $x>_T y$ iff $\alpha>\beta$ and $x\cap \kappa_\beta=y$.

Then $\l T, <_T \r$ is a $\kappa-$tree, since $\kappa$ is an inaccessible.

\begin{lemma}\label{lem1}

$\l T, <_T \r$ has a $\kappa-$branch.

\end{lemma}
\pr
Let $\l x_\gamma \mid \gamma<\kappa\r$ be an enumeration of $T$.
\\There is a club $C\subseteq \kappa$ such that for every $\gamma,\delta\in C$ the following hold:

\begin{enumerate}
  \item $\kappa_\gamma=\gamma$,
  \item the level of $x_\gamma\geq \gamma$,
  \item if $\gamma<\delta$, then the level of $x_\gamma<\delta$.
\end{enumerate}

For every $\gamma\in C$, pick $p_\gamma\in P_\kappa$ such that $p_\gamma\Vdash \lusim{A}\cap \kappa_{\xi_\gamma}=x_\gamma$, where $\xi_\gamma\geq \gamma$ denotes the level of $x_\gamma$.

Now, $\kappa$ is a Mahlo cardinal and an Easton support was used, hence there is a stationary $S\subseteq C$ such that for every $\gamma, \delta\in S$, $p_\gamma$ and $p_\delta$ are compatible.
\\Take any two $\gamma<\delta$ in $S$. Then $x_\delta\cap \xi_\gamma=x_\gamma$ due to the compatibility of $p_\gamma$ and $p_\delta$.
\\So, $\{x_\gamma \mid \gamma\in S\}$ is a $\kappa-$branch.
\\
$\square$

Let now $b=\{x_i \mid i<\kappa\}$ be a maximal $\kappa-$branch in $T$.
For every $i<\kappa$ fix $p_i\in P_\kappa$ which witnesses that $x_i \in T$.
\\By the assumption made, $\bigcup_{i<\kappa} x_i \not = A$, since this union is in $V$.
\\Then for every $i<\kappa$, there is $i'\geq i$ such that $x_{i'}$ is a splitting point of $T$.
\\Denote by $y_{i'}$ an immediate successor of $x_{i'}$ which is not in $b$.
Let $q_{i'}$ be a condition which witnesses that $y_{i'} \in T$.
\\Let $C\subseteq \kappa$ be a club such that for every $i_1,i_2\in C, i_1<i_2$, we have $i_1'<i_2$.

The next lemma provides the desired contradiction, since $P_\kappa$ satisfies $\kappa-$c.c.

\begin{lemma}

The conditions $\{q_{i'}\mid i\in C\}$ are pairwise incompatible.

\end{lemma}
\pr
Let $i_1<i_2$ be in $C$.
Then
$$q_{i_1'} \Vdash \lusim{A}\cap \kappa_{i_1'}\not =x_{i_1'+1}.$$
However,
$$q_{i_2'} \Vdash \lusim{A}\cap \kappa_{i_1'} =x_{i_1'+1},$$
since $y_{i_2'}>_T x_{i_2}>_T x_{i_1'+1}$.
This is possible onle when $q_{i_1'}$ and $q_{i_2'}$ are incompatible.
\\
$\square$

Let us give an   example of the Easton support iteration $P_\kappa$ which adds a fresh subset, however
we give up here the assumption that $|\lusim{Q}_\beta|<\kappa$.

Let $\l \kappa_\beta \mid \beta<\kappa\r$ be an increasing sequence of measurable cardinals above an inaccessible $\kappa$.

Let $\l P_\alpha, \lusim{Q}_\beta\mid \alpha\leq \kappa, \beta<\kappa\r$ be an Easton support iterations of the Prikry  forcings, i.e.
 for each $\beta<\kappa$, $Q_{\beta}$ is the Prikry forcing with a normal ultrafilter over $\kappa_\beta$.

Let $G_\kappa$ be a generic subset of $P_\kappa$.
For every $\beta<\kappa$, let $b_{\beta}$ be the Prikry sequence added by $G_\kappa$ to $\kappa_\beta$.

\begin{lemma}
The set
$$A=\{\alpha<\kappa \mid \text{ the first element of the sequence } b_{\kappa_\alpha} \text{ is } 0\}$$
is a fresh subset of $\kappa$.

\end{lemma}
\pr
Every initial segment of $A$ is in $V$ due to the support used.
On the other hand $A\not \in V$, since every condition in the forcing $P_\kappa$ should be bounded in $\kappa$, and so it can be extended to one which forces $b_{\kappa_\alpha}(0)=0$ or to one forcing $b_{\kappa_\alpha}(0)\not=0$.
\\
$\square$

Let us turn now to a general case, i.e. we assume only that $\kappa$ is an inaccessible.

Our aim  will be to prove the following :

\begin{theorem}\label{thmeastonsupp}
Let $\kappa$ be an inaccessible  cardinal and
$\l P_\alpha, \lusim{Q}_\beta\mid \alpha\leq \kappa, \beta<\kappa\r$ be  the Easton support iteration of Prikry type forcing notions.
Let $I$ be a stationary subset of $ \kappa$ which consists of singular cardinals $\mu$ such that for every $\gamma<\mu, |P_\gamma|<\mu$.
Suppose that for every
 $\alpha\in I$,  $\Vdash_{P_{\alpha}} \langle P_\kappa\setminus {\alpha} , \leq^* \rangle \mbox{ is }\alpha^{++} \mbox{-closed.}$
\\
Let $G_\kappa\subseteq P_\kappa$ be a generic. Then, in $V[G_\kappa]$, there is no fresh subsets of $\kappa$.

\end{theorem}

\begin{remark}
Similar results were proved in  \cite{RestElm} for cardinals above $\kappa$. The proof there is based on the fact that $|P_\kappa|=\kappa$ and it is much easier.

\end{remark}

\pr
Let $G_\kappa$ be a generic subset of $P_\kappa$.

We would like to show that there is no fresh subset of $\kappa$ in $V[G_\kappa]$.

Suppose otherwise. Work in $V$. Let $\lusim{A}$ be a name of such subset and let $\lusim{f}$ be a name of the characteristic function of $A$. Fix some $p\in G_\kappa$ which forces this.

Let $\zeta< \kappa$ be an ordinal for which a new subset is added in the extension from $V$ to $V\left[G\right]$. Such $\zeta $ exists, since the forcings $Q_{\alpha}$ (for $\alpha<\kappa$) have cardinality below $\kappa$, and at least one of them is non-trivial.
By increasing if necessary, we can assume that $\zeta$ is a singular cardinal.
Then  $\l P_\kappa\setminus \zeta, \leq^*\r$  is more than $\zeta$-closed. Thus, there exists a condition $q\in G_\kappa\upr P_{\zeta}$ which forces that a new subset is added to $\zeta$. For simplicity, assume that $p\upr \zeta$  forces this.

Given a condition $r\in P_\kappa$, let us denote by $r({\gamma})$ its $\gamma-$th coordinate,\\ i.e. $r=\l r(\gamma)\mid \gamma<\kappa\r$.

We divide into two cases as in Theorem \ref{thm4-1}.

\textbf{Case 1.} There exists $ \mu\in \left( \zeta, \kappa \right)\cap I$ and a condition $ p^* \in P_{\mu} $ which forces that the following property holds:
\begin{align*}
	&\exists p\in G_\mu \  \exists s\in P\setminus \mu \forall r\geq^* s \  \exists \xi <\kappa \  \exists r_0,r_1\geq^* r  , \\
	& V\vDash \left( {p}^{\frown} r_0 \parallel {\lusim{f}}\restriction { \xi }  , \  {p}^{\frown} r_1 \parallel {\lusim{f}}\restriction { \xi } \right) \mbox{, and the decisions are different.}
\end{align*}

\textbf{Case 2.} For every $ \mu\in \left( \zeta, \kappa \right)\cap I $, every condition in $ P_{\mu} $ forces that--
\begin{align*}
	& \forall p\in \lusim{G}_\mu \ \forall s\in P_\kappa\setminus \mu \exists r\geq^* s \  \forall \xi<\kappa  \  \forall r_0,r_1\geq^* r, \\
	& V\vDash \mbox{ If }  {p}^{\frown} r_0 \parallel {\lusim{f}}\restriction { \xi }  \mbox{ and } \  {p}^{\frown} r_1 \parallel {\lusim{f}}\restriction { \xi }  \mbox{ then the decisions are the same.}
\end{align*}

The treatment of the first case is exactly as in Theorem \ref{thm4-1}.
Let us deal with the second case.

As in \ref{thm4-1}, we define for every  $ \mu<\kappa$ 
\begin{align*}
	e\left(\mu\right) = \{ & r\in P_\kappa\setminus \mu \colon \forall p\in G_{\mu} \     \forall \xi<\kappa \  \forall r_0,r_1\geq^* r, \\
	& V\vDash \mbox{ If }  {p}^{\frown} r_0 \parallel {\lusim{f}}\restriction{ \xi }  \mbox{ and } \  {p}^{\frown} r_1 \parallel {\lusim{f}}\restriction{ \xi }  \mbox{ then the decisions are the same.}   \}
\end{align*}

By Claim 1 of \ref{thm4-1}, it is $\leq^*-$dense open subset of $P_\kappa\setminus \mu$, for every $\mu\in (\zeta, \kappa)\cap I$.
Again, here is the only place where $\mu^{++}-$closure of the direct order on $P_\kappa\setminus \mu$ is used.

Given a generic $G_\kappa\subseteq P_\kappa, p\in G_\kappa$, define
 in $V[G_\kappa]$, $$S=\{\xi<\kappa \mid  \xi \text{ ia limit ordinal and } \exists t\geq^*p,t\in G_\kappa \text{ such that }
t\upr \xi \Vdash_{P_\xi} t\setminus \xi||_{ P_\kappa\setminus{\xi}}
\lusim{A}\cap \xi)\}.$$

For every $\xi\in S$, fix some $t^\xi\in G_\kappa$ such that  $t^\xi \geq^*p$ and
$t^{\xi}\upr \xi \Vdash_{P_\xi} t^\xi\setminus \xi||_{ P_\kappa\setminus{\xi}} \lusim{A}\cap \xi$.
So,   there is $p^\xi \in G\upr P_\xi,
p^\xi\geq t^\xi\upr \xi, (p^\xi)^\frown t^\xi\setminus \xi\Vdash \lusim{A}\cap \xi=a_\xi$, for some $a_\xi\in V$.
Then there is a finite $b^\xi\subseteq \xi$ such that $p^\xi\setminus b^\xi \geq^* p\upr (\xi\setminus b^\xi)$.


Suppose for a moment that $S$ is stationary in $V[G_\kappa]$.
\\
Then we can find a stationary subset $S'$ of $S$ and a finite $b$ such that for every $\xi\in S'$, $b^\xi=b$.
\\Now we can freeze $p^\xi\upr \max(b)$.
Denote $\max(b)$ by $\mu^*$.

Let $\mu\in I$ be a cardinal above $\mu^*$.
\\Consider the  set  $e(\mu)$ defined above.
It is $\leq^*-$dense open subset of $P_\kappa\setminus \mu$,  in $V[G_\mu]$ above $p\setminus \mu$.
In particular there is $r\in P_\kappa\setminus \mu, r\geq^* p\setminus \mu$ such that
\\for every $p'\in G_\mu$, for every $\xi<\kappa$ and for every $r_0,r_1\geq^* r$
$$V\models \text{ If } p'{}^\frown r_0 || \lusim{A} \cap \xi \text{ and }  p'{}^\frown r_1 || \lusim{A} \cap \xi \text{ then the decisions are the same}.$$
Recall that $t^\xi\geq^* p$ and $p^\xi\setminus \mu^* \geq^* p\upr (\mu^*, \xi)$, for every $\xi\in S'$.
If we were able to conclude from this that $p^\xi\setminus \mu \geq^* r$, then it will imply that  $A\in V$.
However it need not be the case since the support of $r$ may be bigger than those of $p\setminus \mu$ and incompatibility may occur on coordinates outside of $\supp(p)$.

Let us argue that it is possible to overcome this obstacle.
\\
Work in $V$. Set $p_0=p$. Use Lemma \ref{Lem2-3} to
find $q\geq^* p$ and $\mu_0\geq \sup(\supp(p)), \mu_0\in I$ such that $q\Vdash \lusim{\mu}^*\leq \mu_0$.
\\Then, we use the density of $e(\mu_0)$ to find $r_0\geq^* q\setminus \mu_0$ such that
$$q\upr \mu_0\Vdash_{P_{\mu_0}} r_0\in e(\mu_0).$$
Set $p_1=q\upr \mu_0{}^\frown r_0$.
\\Next, we run the argument above with $p=p_0$ replaced by $p_1$.
Again, using Lemma \ref{Lem2-3}
find $q_1\geq^* p_1$ and $\mu_1\geq \sup(\supp(p_1)),\mu_0+1$,$\mu_1\in I$  such that $q_1\Vdash \lusim{\mu}^*\leq \mu_1$, where $\mu^*$ is now defined using $p_1$ instead of $p$.
\\Then,  use the density of $e(\mu_1)$ to find $r_1\geq^* q_1\setminus \mu_1$ such that
$$q_1\upr \mu_1\Vdash_{P_{\mu_1}} r_1\in e(\mu_1).$$
Set $p_2=q\upr \mu_1{}^\frown r_1$.
\\Continue by induction and define $p_n,q_n, r_n, \mu_n$, for every $n<\omega$.
\\Finally set $\mu_\omega=\bigcup_{n<\omega} \mu_n$ and $p_\omega=\bigcup_{n<\omega} p_n = q_\omega=\bigcup_{n<\omega} q_n $.
Then, for every $n<\omega$, $$p_\omega\upr \mu_n \Vdash p_\omega\setminus \mu_n \in e(\mu_n),$$
since $p_\omega\setminus \mu_n \geq^* r_n\in e(\mu_n)$ and $e(\mu_n)$ is dense open. Also, $\sup(\supp(p_\omega))=\mu_\omega$.
\\Pick now a generic $G\subseteq P_\kappa$ with $p_\omega\in G$. Let $S_{p_\omega}$ be defined as $S$ above only with $p_\omega$ replacing $p$.
Assuming its stationarity, define  $q^\xi\geq^* p_\omega , p^\xi\in G\upr P_\xi$ for $\xi\in S'$ exactly as above.
Then there will be a stationary $S'\subseteq S$ and $n^*<\omega$ such that $p^\xi\upr (\mu_{n^*}, \xi)\geq^* p_\omega\upr (\mu_{n^*}, \xi)$, since
$\sup(\supp(p_\omega))=\mu_\omega$ and a non-direct extension is used at finitely many places only.
\\
Shrink $S'$ further to $S''$ and stabilise the value of the function $\xi \mapsto   p^\xi\upr \mu_{n^*}^{+}$.
\\Finally, we use that $p_\omega\setminus \mu_{n^*} \in e(\mu_{n^*})$.

Hence, the following lemma will complete the proof.
We prove it for the initial $p$, but the same argument works for $S_{p_\omega}$ or any $S_x$ with $x\geq p$.

\begin{lemma}\label{lem1-1}
	
	$S$ is stationary in $V[G_\kappa]$.

\end{lemma}
\pr
The argument will be similar to those of \ref{thm101}.
\\Suppose otherwise.
Let $C\subseteq \kappa$ be a club disjoint from $S$.
Assume that $p\in G_\kappa$  forces this, otherwise replace it by a stronger condition doing this.
\\Work in $V$. Pick an elementary submodel $M\preceq H_\chi$ such that

\begin{enumerate}
	\item $|M|=\delta<\kappa$,
	\item $M\cap \kappa=\delta$,
\item $\cof(\delta)<\delta$,
	\item ${}^{\cof(\delta)>}M\subseteq M$,
	\item $\kappa,P_\kappa,  \lusim{C}, p\in M$.
\end{enumerate}

Pick a cofinal in $\delta$ sequence $\l \delta_i \mid i<\cof(\delta) \r$ consisting of singulars and with $\delta_0>\cof(\delta)$.

Let $G_{\cof(\delta)+1}\subseteq P_{\cof(\delta)+1}$ be a generic with $p\upr {\cof(\delta)+1}\in G_{\cof(\delta)+1}$.
\\Work in $V[G_{\cof(\delta)+1}]$.

Consider
$$D^0=\{r\geq p\mid \exists c<\kappa \quad r\Vdash c=\min(\lusim{C}\setminus \delta_0)\}.$$

Clearly, $D^0$ is a dense open and it belongs to $M$.
\\
Apply Lemma \ref{lem101}. Then there will be $\alpha_0<\kappa$ and $q_0\geq^* p, q\upr {\cof(\delta)+1}=p\upr {\cof(\delta)+1}$, in $M$ such that for every $r\in P_{\alpha_0}, r\geq_{P_{\alpha_0}} q_0\upr \alpha_0$ there is $r'\geq_{P_{\alpha_0}} r$ such that
$r'{}^\frown q_0\setminus \alpha_0 \in D_0$.
So, for every such $r'$ there is $c(r')<\kappa$ such that
$$r'{}^\frown q_0\setminus \alpha_0 \Vdash c(r')=\min(\lusim{C}\setminus \delta_0).$$
Note that $P_{\alpha_0}\subseteq M$. Hence all $c(r')$'s are in $M$.
Also $|P_{\alpha_0}|<\kappa$. Hence their sup is below $\kappa$, and then, by elementarity, in $M$.
Denote it by $c_0^*$.

Define, for every $\tau<\kappa$,

$$D_\tau=\{r \in P \mid r\geq^* p \text{ and } \exists \xi\in \kappa\setminus \tau \quad 
r\upr \xi \Vdash_{P_\xi} (r\setminus \xi)||_{ P_{>\xi}} \lusim{A}\cap \xi)\}.$$

\begin{claim}
$D_\tau$ is $\leq^*-$dence open above $p$.

\end{claim}
\pr
Set $\tau_0=\tau+1$. Consider

$$D(\tau_0)=\{r\in P_\kappa \mid r\perp p \text{ or } (r\geq p \text{ and }  r\parallel \lusim{A} \cap \tau_0)\}.$$

By Lemma \ref{lem101}, there are $\alpha_0$ and $q\geq^* p$ such that

$$q\upr \alpha_0 \Vdash \exists b\in \lusim{G}_{\alpha_0} \quad b^\frown q\setminus \alpha_0 \parallel \lusim{A} \cap \tau_0.$$

If $\alpha_0\leq \tau_0$ then we are done. Suppose that $\alpha_0> \tau_0$. Consider $D(\alpha_0)$ and again, using Lemma \ref{lem101}, pick $\alpha_1$ and $q_1\geq^* q$
such that

$$q_1\upr \alpha_1 \Vdash \exists b\in \lusim{G}_{\alpha_1} \quad b^\frown q_1\setminus \alpha_1 \parallel \lusim{A} \cap \alpha_0.$$

If $\alpha_1\leq \alpha_0$, then $q_1\in D(\alpha_0)$ and we are done. If $\alpha_1> \alpha_0$, then continue and define in the same fashion $\alpha_2, q_2$ etc.
\\ Suppose that the process continues infinitely many steps. Then we will have
$$\alpha_0<\alpha_1<...<\alpha_i<... \text{ and } q_0\leq^* q_1\leq^*...\leq^*q_i\leq^*..., i<\omega.$$
Let $\alpha^*=\bigcup_{i<\omega} \alpha_i$ and $q^*\geq^* q_i$, for every $i<\omega$.
Let $G\subseteq P_\kappa$ be a generic with $q^*\in G$. For every $i<\omega$, let $G_{\alpha_i}=G\upr P_{\alpha_i}$ and $G^*=G\upr P_{\alpha^*}$.
\\Now, for every $i<\omega$, there are $r_{i+1}\in G_{\alpha_{i+1}}$ and $a_i\in V$ such that
$$r_{i+1}{}^\frown q_{i+1}\setminus \alpha_{i+1}\Vdash \lusim{A}\cap \alpha_i=a_i.$$

Then,

$$r_{i+1}{}^\frown q^*\setminus \alpha_{i+1}\Vdash \lusim{A}\cap \alpha_i=a_i.$$

Set $a=\bigcup_{i<\omega}a_i$.
Then, in $V[G_{\alpha^*}]$, $q^*\setminus \alpha^*\Vdash \lusim{A}\cap \alpha^*=a$, since of $r_i$'s is in $G_{\alpha^*}$.
So, there is $r \in G_{\alpha^*}$ such that
$$r^\frown q^*\setminus \alpha^*\Vdash  \lusim{A}\cap \alpha^*=a.$$
Hence,
$$r\Vdash_{P_{\alpha^*}}(q^*\setminus \alpha^*||\lusim{A}\cap \alpha^*).$$

The only requirement  on $G_{\alpha^*}$ is that $q^*\upr \alpha^*\in G_{\alpha^*}$. Hence,

$$q^*\upr \alpha^* \Vdash q^*\setminus \alpha^*|| \lusim{A}\cap \alpha^*.$$
So, $q^* \in D_{\tau}$.
\\
$\square$ of the claim.

Consider $D_{c_0^*}$.

 It is in $M$, as well.
So, inside $M$, we can pick $\xi_0\geq c_0^*,\xi_0\in I $ and $t_0\geq^* p $ such that
$$
t_0\upr \xi_0 \Vdash_{P_{\xi_0}}  t_0\setminus \xi_0||_{ P_\kappa\setminus {{\xi_0}}} \lusim{A}\cap \xi_0)\}.$$

We continue the same process only with $c_0$ replaced by $\max(\xi_0, \delta_1)$ and $p$ by $q_0$.

At limit stages of the construction and at the stage $\cof(\delta)$ itself we would like to put $\l q_i\mid i<\cof(\delta)\r$ into a single condition $q^*$.
This is possible since $\delta_0>\cof(\delta)$, and so,  we have enough completeness to proceed.

We have  $p\upr {\cof(\delta)+1}^\frown q^*\Vdash \delta\in \lusim{C}$, since $C$ is forced to be closed by $p$.

Let $G\subseteq P_\kappa$ be a generic with $q^*\in G$.
Set $G_\delta=G\upr P_\delta$.
Then for every $  i<\cof(\delta)$, $q_i\upr \delta\in G_\delta$.
Hence, there are $r_i\in G_{\xi_i}, r_i\geq q_i\upr \xi_i$ and $a_i\in V$ such that
$$r_i{}^\frown q_i\setminus \xi_i \Vdash \lusim{A}\cap \xi_i=a_i.$$
Set $a=\bigcup_{i<\cof(\delta)}a_i$. Then, in V[G], $A\cap \delta=a$, since each $r_i{}^\frown q_i\setminus \xi_i$ is in $G$.
\\Remember that the only requirement on a generic set $G$ was that $q^*$
belongs to it. Hence, back in $V$,

$$q^*\upr \delta \Vdash_{P_\delta} q^*\setminus \delta||_{ P_{>\delta}} \lusim{A}\cap \delta.$$

So, $\delta\in S$ and as it was shown above $\delta\in C$ as well. Contradiction.
\\
$\square$

\section{ No fresh subsets of $\kappa$ in the full support }

Our aim  will be to prove the following :

\begin{theorem}
Let $\kappa$ be an inaccessible  cardinal and
$\l P_\alpha, \lusim{Q}_\beta\mid \alpha\leq \kappa, \beta<\kappa\r$ be the full iteration of Prikry type forcing notions.
Suppose that for every $\beta<\kappa$, for every $x,y,z\in \lusim{Q}_\beta$, if $z\leq_{\lusim{Q}_\beta}^*x,y$ and $x,y$ are compatible according to $\leq_{\lusim{Q}_\beta}$, then they are compatible according to $\leq_{\lusim{Q}_\beta}^*$, i.e. there is $e \in \lusim{Q}_\beta, e\geq_{\lusim{Q}_\beta}^*x,y$.\footnote{ Note that if $\leq=\leq^*$, then this holds trivially. Prikry, Magidor, Radin forcings, their supercompact versions, etc., have this property. Actually, any reasonable forcing of this type has this property.  }

Let $G_\kappa\subseteq P_\kappa$ be a generic. Then, in $V[G_\kappa]$, there is no fresh subsets of $\kappa$.

\end{theorem}

\begin{remark}
More restrictive results were proved in  \cite{RestElm} for cardinals above $\kappa$.
The present proof can be easily modified for higher cardinals.

\end{remark}

\pr
Let $G_\kappa$ be a generic subset of $P_\kappa$.

We would like to show that there is no fresh subset of $\kappa$ in $V[G_\kappa]$.

Suppose otherwise. Work in $V$. Let $\lusim{A}$ be a name of such subset and let $\lusim{f}$ be a name of the characteristic function of $A$. Fix some $p\in G_\kappa$ which forces this.



Let $\zeta< \kappa$ be an ordinal for which a new subset is added in the extension from $V$ to $V\left[G\right]$. Such $\zeta $ exists, since the forcings $Q_{\alpha}$ (for $\alpha<\kappa$) have cardinality below $\kappa$, and at least one of them is non-trivial.
By increasing if necessary, we can assume that $\zeta$ is a singular cardinal.
Then  $\l P_\kappa\setminus \zeta, \leq^*\r$  is more than $\zeta$-closed. Thus, there exists a condition $q\in G_\kappa\upr P_{\zeta}$ which forces that a new subset is added to $\zeta$. For simplicity, assume that $p\upr \zeta$  forces this.

Let $I$ be a subset of $\kappa$ which consists of singular cardinals $\tau$ such that for every $\rho<\tau, |P_\rho|<\tau$.


As before,  we divide into two cases:

\textbf{Case 1.}  There exists  $ \mu\in \left( \zeta, \kappa \right)\cap I$ and a condition $ p^* \in P_{\mu}, p^*\geq p\upr \mu $ which forces that the following property holds:
\begin{align*}
	&\exists p'\in \lusim{G}_\mu \  \exists s\in P_\kappa\setminus \mu, s\geq p\setminus \mu \forall r\geq^* s,  r({\mu})=s({\mu}), r({\mu^+})=s({\mu^+}) \  \exists \xi <\kappa \  \exists r_0,r_1\geq^* r,
\\
& r_0({\mu})=r_1({\mu})=r({\mu}),
r_0(\mu^+)=r_1({\mu^+})=r({\mu^+}),  \\
	& V\vDash \left( {p'}^{\frown} r_0 \parallel {\lusim{f}}\restriction { \xi }  , \  {p'}^{\frown} r_1 \parallel {\lusim{f}}\restriction { \xi } \right) \mbox{, and the decisions are different,}
\end{align*}
where $ G_{\mu} $ denotes the canonical name for the generic set for $ P_{\mu} $. \\By extending $ p^* $, if necessary, we can decide the value of $ p' $ in the statement above, and thus assume that $ p^* \geq p' $. Let $ \lusim{s} $ be a $ P_{\mu} $-name for $ s $ from the above property, and assume that this is forced by $ p^* $.

Note that here we do not assume $\mu^{++}-$completence of the direct order. Additional requirements  are included in order to compensate this.
Still, the treatment of this case repeats completely Case 1 of \ref{thm4-1}.

\textbf{Case 2.} For every $ \mu\in \left( \zeta, \kappa \right)\cap I$, every condition in $ P_{\mu} $ stronger than $p\upr \mu$ forces ( and so, $p\upr \mu$ forces) that--
\begin{align*}
	& \forall p'\in \lusim{G}_\mu \ \forall s\in P\setminus \mu, s\geq p\setminus \mu \exists r\geq^* s, r({\mu})=s({\mu}), r({\mu^+})=s({\mu^+}) \  \forall \xi<\kappa  \  \forall r_0,r_1\geq^* r, \\
& r_0({\mu})=r_1({\mu})=r({\mu}),
r_0({\mu^+})=r_1({\mu^+})=r({\mu^+}), \\
	& V\vDash \mbox{ If }  {p'}^{\frown} r_0 \parallel {\lusim{A}}\cap { \xi }  \mbox{ and } \  {p'}^{\frown} r_1 \parallel {\lusim{A}}\cap { \xi }  \mbox{ then the decisions are the same.}
\end{align*}
For every  $ \mu\in \left( \zeta, \kappa \right)\cap I$, define (in $V[G_\mu]$)
\begin{align*}
	e'\left(\mu\right) = \{ & r\in P_\kappa\setminus \mu \colon \forall p'\in G_{\mu} \     \forall \xi<\kappa \  \forall r_0,r_1\geq^* r, \\
& r_0({\mu})=r_1({\mu})=r({\mu}),
r_0({\mu^+})=r_1({\mu^+})=r({\mu^+}), \\
	& V\vDash \mbox{ If }  {p'}^{\frown} r_0 \parallel {\lusim{f}}\restriction{ \xi }  \mbox{ and } \  {p'}^{\frown} r_1 \parallel {\lusim{f}}\restriction{ \xi }  \mbox{ then the decisions are the same.}   \}
\end{align*}

\begin{claim} For every $ \mu\in \left( \zeta, \kappa \right)\cap I$,

$e'(\mu)$ is $ \leq^* $-dense  in $ P\setminus \mu$ (above $p\setminus \mu$).\footnote{ It need not be open.}

\end{claim}
\pr
Note that  $ \left| G_{\mu} \right|\leq \mu^+$, the forcing $ \l P_\kappa\setminus \mu^{++} \r$ is $\mu^{++}-$closed and the coordinates $\mu, \mu^+$ do not change.
Hence, for every given $r'\in  P_\kappa\setminus \mu$, we can construct a $\leq^*-$increasing sequence of a length $|G_\mu|$ of conditions stronger than $r'$,
which takes care of each $p'\in G_\mu$. Then its upper bound will be in $e(\mu)$.
\\
$\square$ of the claim.

Given a generic $G_\kappa\subseteq P_\kappa, p\in G_\kappa$, define
 in $V[G_\kappa]$, $$S=\{\xi<\kappa \mid  \xi \text{ ia limit ordinal and } \exists t\geq^*p,t\in G_\kappa \text{ such that }
t\upr \xi \Vdash_{P_\xi} t\setminus \xi||_{ P_\kappa\setminus{\xi}}
\lusim{A}\cap \xi)\}.$$

For every $\xi\in S$, fix some $t^\xi\in G_\kappa$ such that  $t^\xi \geq^*p$ and
$t^{\xi}\upr \xi \Vdash_{P_\xi} t^\xi\setminus \xi||_{ P_\kappa\setminus{\xi}} \lusim{A}\cap \xi$.
So,   there is $p^\xi \in G\upr P_\xi,
p^\xi\geq t^\xi\upr \xi, (p^\xi)^\frown t^\xi\setminus \xi\Vdash \lusim{A}\cap \xi=a_\xi$, for some $a_\xi\in V$.
Then there is a finite $b^\xi\subseteq \xi$ such that $p^\xi\setminus b^\xi \geq^* p\upr (\xi\setminus b^\xi)$.


The argument of Lemma \ref{lem1-1} applies without changes in the present situation, and shows
 that $S$ is stationary (in $V[G_\kappa]$).

Then we can find a stationary subset $S'$ of $S$ and a finite $b$ such that for every $\xi\in S'$, $b^\xi=b$.
\\Now we can freeze $p^\xi\upr \max(b)^{++}$.
\\If the following holds,

$(\aleph)$: \emph{For every $\beta<\kappa$, for every $s,t,r \in \lusim{Q}_\beta$, if $s\leq_{\lusim{Q}_\beta}^*t,r$, then there is $e \in \lusim{Q}_\beta, e\geq_{\lusim{Q}_\beta}^*t,r$}.\footnote{ For example if $P_\kappa$ is the Magidor iteration of Prikry forcings.
Also note that we do not to split into Cases 1,2 in this type of situation. }

then it is easy  to conclude that $A\in V$, which is impossible.
\\Deal with the general situation.
Let $\mu\in I$ be a cardinal $\geq \max(b)$.
Freeze $p^\xi\upr \mu^{++}+1$.
\\Then for every $\xi,\rho \in S'$,
$$p^\xi(\mu)=p^\rho(\mu) \text{ and } p^\xi(\mu^+)=p^\rho(\mu^+).$$

There is $s\in G_\kappa, s\geq p, s\upr \mu^{++}\geq p^\xi\upr\mu^{++}, \xi\in S'$,  such that
$$s\Vdash \forall \xi\in \lusim{S}' \exists t^\xi\geq^* p, t^\xi\in \lusim{G}_\kappa \quad\exists p^\xi \in \lusim{G}_\kappa\upr P_\xi \quad b_\xi=b \text{ which has the properties  above}.$$
\\Consider the  set  $e'(\mu)$ defined above.
By the claim,
it is $\leq^*-$dense  subset of $P_\kappa\setminus \mu$,  in $V[G_\mu]$ above $p\setminus \mu$.
In particular there is $r\in P_\kappa\setminus \mu, r\geq^* s\setminus \mu$ such that
\\for every $p'\in G_\mu$, for every $\xi<\kappa$ and for every $r_0,r_1\geq^* r$
$$V\models \text{ If } p'{}^\frown r_0 || \lusim{A} \cap \xi \text{ and }  p'{}^\frown r_1 || \lusim{A} \cap \xi \text{ then the decisions are the same}.$$
Pick now a generic $G'\subseteq P_\kappa$ with $G_\mu=G'\upr P_\mu$ and $ s\upr \mu^\frown r\in G'$.

We assumed the following,

$(\beth)$: \emph{For every $\beta<\kappa$, for every $x,y,z\in \lusim{Q}_\beta$, if $z\leq_{\lusim{Q}_\beta}^*x,y$ and $x,y$ are compatible according to $\leq_{\lusim{Q}_\beta}$, then they are compatible according to $\leq_{\lusim{Q}_\beta^*}$, i.e. there is $e \in \lusim{Q}_\beta, e\geq_{\lusim{Q}_\beta}^*x,y$.}

Now,  $P_\kappa$ is a full support iteration. So,
$(\beth)$ implies -

$(\gimel)$: \emph{For every $\beta<\kappa$, for every $u,v,w\in P_\kappa$, if $w\leq_{P_\kappa}^*u,v$ and $u,v$ are compatible according to $\leq_{P_\kappa}$, then they are compatible according to $\leq_{P_\kappa}^*$.\footnote{ It need not be the case for Easton or non-stationary support, since then the support of $w$ may be strictly smaller than those of $u$, $v$, and $u,v$ may disagree on a common coordinate outside the support of $w$.}
}

It follows, in this case,  that $A\in V$. Contradiction.
\\
$\square$


\newpage
\begin{thebibliography}{AA}

\bibitem{ben2017homogeneous} Omer Ben-Neria and Spencer Unger, Homogeneous changes in cofinality with applications to HOD, Journal Math. logic, vol 17 (2017) no.2

\bibitem{friedman2009number} Sy Freidman and Menachem Magidor, The JSL, vol 74, (2009), pp. 1069-1080.

\bibitem {hb} Moti Gitik.
\newblock {P}rikry-type forcings.
\newblock In {\em Handbook of set theory}, pages 1351--1447. Springer, 2010.

\bibitem{RestElm}  Moti Gitik and Eyal Kaplan.
\newblock On restrictions of ultrafilters from generic extensions to ground
models.
\newblock {\em The Journal of Symbolic Logic}, pages 1--31, 2021.

\bibitem{G-K-nonst} Moti Gitik and Eyal Kaplan.
\newblock Non-stationary support iterations of prikry forcings and restrictions
of ultrapower embeddings to the ground model.
\newblock {\em Annals of Pure and Applied Logic}, 174(1):103164, 2023.


\bibitem{G-K-easton} Moti Gitik and Eyal Kaplan.
\newblock On {E}aston support iteration of {P}rikry-type forcing notions.
\newblock {\em arXiv preprint arXiv:2301.12421}, 2023.

\bibitem{hamkins2001gap} Joel~David Hamkins.
\newblock Gap forcing.
\newblock {\em Israel Journal of Mathematics}, 125(1):237--252, 2001.


\bibitem{K-full} Eyal Kaplan.
\newblock The {M}agidor iteration and restrictions of ultrapowers to the ground
model.
\newblock {\em arXiv preprint arXiv:2202.04980}, 2022.



\bibitem{Eyal-thesis} Eyal Kaplan, PhD Thesis, Tel Aviv University, 2023.

\bibitem{Mag} Menachem Magidor, How large is the first strongly compact cardinal? Or study of identity crises, Ann. Math. Logic, vol. 10 (1976), pp. 33-57.

\bibitem{Levy-Solovay} Azriel Levy and Robert M. Solovay. Measurable cardinals and the continuum hypothesis. Israel Journal of Mathematics 5:234-248, 1967.


\end{thebibliography}
\end{document}